\documentclass[11pt]{article}

\usepackage[english]{babel}
\usepackage{latexsym}
\def\nl2#1{\left\| #1\right\|_{L^2(\O)}}

\def\bel{\begin{equation}\label}
\def\beq{\begin{equation}}
\def\ee{\end{equation}}
\usepackage{fourier}
\usepackage{bbm}
\usepackage[pdftex]{graphicx}
\usepackage{amsfonts}
\usepackage{amsmath}
\usepackage{amsthm}
\usepackage{amssymb}
\usepackage{empheq}
\usepackage{textcomp}
\usepackage{color}
\usepackage{caption}
\usepackage{tikz}
\usepackage{mathrsfs}
\usepackage{graphicx}
\usepackage[onehalfspacing]{setspace}
\usepackage[toc,page]{appendix}
\usepackage{cases}

\DeclareGraphicsExtensions{.png}

\def\O{\Omega}

\newcommand{\bela}{\begin{equation} \label}
\newcommand{\eeq}{\end{equation}}
\newcommand{\ba}{\begin{array}}
\newcommand{\ea}{\end{array}}

\makeatother

\newtheorem{theorem}{Theorem}[section]
\newtheorem{proposition}{Proposition}[section]
\newtheorem{lemma}{Lemma}[section]
\newtheorem{corollary}{Corollary}[section]
\newtheorem{remark}{Remark}[section]
\newtheorem{notation}{Notation}[section]

\catcode`\@=12

\begin{document}


\begin{center}
{\Large \bf  Spectral Asymptotic for the Infinite Mass Dirac Operator in bounded domain}
\end{center}
\begin{center}

Badreddine Benhellal\footnote{\textit{$2010$ Mathematics Subject Classification}.	81Q10 ,	81V05, 35P15, 58C40}{\*}{\*}
\end{center}

\medskip

\noindent

\smallskip

{\small
\begin{abstract}
In this paper, we study  a singular perturbation of a problem used in dimension two to model graphene or in dimension three to describe the quark confinement phenomenon in hadrons. The operators we consider are of the form $\textbf{H} +M\beta V(x)$,  where $\textbf{H}$ is the free Dirac operator, $\beta$ is a constant matrix, $V(x)$ is a real valued piecewise constant potential having a jump discontinuity across a smooth interface and $M$ is the mass that we can see as a coupling constant. In particular, we perform a complete  asymptotic expansion  of spectral quantities as the mass $M$ tends to $+\infty$.
\end{abstract}
}

\tableofcontents

\section{Introduction} \label{sec1}

 Let $\Omega$ be a smooth bounded connected domain of $\mathbb{R}^n$ with $n\in \lbrace 2,3\rbrace$. We define $T_2$ and $T_3$ the differential expressions associated to the massless Dirac operator on $\mathbb{R}^2$ and $\mathbb{R}^3$, respectively as
\begin{align*}
&T_2 :=\frac{1}{\textit{i}} \boldsymbol{\sigma}\cdot\nabla_2:= \frac{1}{\textit{i}}\sum_{j=1}^{2}\sigma_{j}\partial_{j}=-\textit{i}\begin{pmatrix}
0 & \partial_1 -\textit{i} \partial_2\\
\partial_1 +\textit{i}\partial_2 & 0
\end{pmatrix}=\begin{pmatrix}
0 & D_2^\ast\\
D_2 & 0
\end{pmatrix},\\
&T_3 :=\frac{1}{\textit{i}} \boldsymbol{\alpha}\cdot\nabla_3:=\frac{1}{\textit{i}} \sum_{j=1}^{3}\alpha_{j}\partial_{j}=\begin{pmatrix}
0 & D_3^\ast\\
D_3 & 0
\end{pmatrix},
\end{align*}
where $\nabla_n$ is the gradiant operator in $\mathbb{R}^n$, the Dirac matrices $\boldsymbol{\alpha} = (\alpha_1,\alpha_2,\alpha_3)$ and $\beta$ are $4\times4$ Hermitian and unitary, given by
$$
\beta=\begin{pmatrix}
1_2 & 0\\
0 & -1_2
\end{pmatrix},\quad
\alpha_k=\begin{pmatrix}
0 & \sigma_k\\
\sigma_k & 0
\end{pmatrix},\quad \text{ for } k=1,2,3.
 $$
Here $\boldsymbol{\sigma} = (\sigma_1, \sigma_2)$ and the Pauli matrices $\sigma_1$, $\sigma_2$ and $\sigma_3$ are defined by
$$ \sigma_1=\begin{pmatrix}
0 & 1\\
1 & 0
\end{pmatrix},\quad \sigma_2=
\begin{pmatrix}
0 & -i\\
i & 0
\end{pmatrix} ,\\
\quad
\sigma_3=\begin{pmatrix}
1 & 0\\
0 & -1
\end{pmatrix}.$$ 

 The free Dirac operator $(\textbf{H},\text{Dom}(\textbf{H}))$ associated with energy of relativistic particle of mass $M>0$ and spin $1/2$, acting on $L^2(\mathbb{R}^3;\mathbb{C}^4)$ in three dimension case (on $L^2(\mathbb{R}^2;\mathbb{C}^2)$ for the 2D case) is defined in the sense of distributions by 
 \begin{itemize}
 \item  in 2D : $ \textbf{H} := T_2 + M \sigma_3$, with Dom($\textbf{H}$) =$H^1(\mathbb{R}^2)$:= $H^1(\mathbb{R}^2;\mathbb{C}^2)$ and $D_2^\ast$ is the adjoint of $D_2$.
 \item  in 3D : $ \textbf{H} := T_3 + M \beta$, with Dom($\textbf{H}$) = $H^1(\mathbb{R}^3)$:=$H^1(\mathbb{R}^3;\mathbb{C}^4)$ and $D_3 = D_3^\ast$ (see appendix for the explicit formula).
 \end{itemize}
 Here $H^1(\mathbb{R}^n)$ is the usuel Sobolev space. For the sequel we use the notations $\Omega_+ :=\Omega$, $\Omega_- :=\mathbb{R}^n\setminus\overline{\Omega}$, $\partial\Omega :=\partial\Omega_\pm$, $\Theta_2 := \sigma_3$, $\Theta_3 :=\beta $, $\boldsymbol{\Lambda}_2:=\boldsymbol{\sigma}$ and $\boldsymbol{\Lambda}_3 := \boldsymbol{\alpha}$ according to the dimension. Then we define the self-adjoint Dirac operator $(\textbf{H}_M,\text{Dom}(\textbf{H}_M))$ with high scalar potential barrier in $n$-dimension by 
 \begin{align}
 \textbf{H}_M := T_n + M\Theta_n(1-\mathbbmss{1}_{\Omega_+}), \text{ with }\text{Dom}(\textbf{H}_M)= H^1(\mathbb{R}^n),
 \end{align} 
where $\mathbbmss{1}_{\Omega_+}$ is the characteristic function of $\Omega_+$. Thus one can see this operator as a compactly supported perturbation of the free Dirac operator. As a consequence,  $\textbf{H}_M$ has purely discrete spectrum between $(-M,M)$ and its essential spectrum is equal to $(-\infty,-M]\cup [M,+\infty)$ (see \cite{STSV},\cite{ATR}). In order to describe the boundary conditions of the limit operator $\textbf{H}_\infty$, we denote by $\textbf{n}_{\Omega_+}: \partial\Omega\rightarrow\mathbb{S}^{n-1}$ the outward normal vector of $\Omega_+$.\/

 Let $\mathcal{D}^n_\infty = \lbrace \psi \in H^1(\Omega_+) : \mathcal{B}^n_{\Omega_+}\psi =\psi \text{ on } \partial\Omega\rbrace,$ with $\mathcal{B}^n_{\Omega_+}(x) = -i\Theta_n\left(\boldsymbol{\Lambda}_n\cdot\textbf{n}_{\Omega_+}(x)\right)$, for all $x\in\partial\Omega$.
 Then the limit operator is defined by
\begin{align*}
 \textbf{H}_\infty :\mathcal{D}^n_\infty&\subset L^2(\Omega_+; \mathbb{C}^n)\rightarrow  L^2(\Omega_+; \mathbb{C}^n),\\
&\textbf{H}_\infty\psi =T_n\psi .
\end{align*}
To identify the boundary conditions, let's take $\psi\in \mathcal{C}^\infty(\overline{\Omega}_+;\mathbb{C}^n)$, then using the Green formula and the definition of the matrix $\mathcal{B}^n_{\Omega_+}$, we get 
\begin{align*}
\langle \psi, T_n\psi\rangle =\langle T_n\psi, \psi\rangle -\langle \psi, \Theta_n\mathcal{B}^n_{\Omega_+} \psi\rangle_{\partial\Omega}.
\end{align*}
 Thus any self-adjoint realization of $T_n$ in $L^2(\Omega_+,\mathbb{C}^n)$ must satisfy: $ \langle \psi, \Theta_n\mathcal{B}^n_{\Omega_+} \psi\rangle_{\partial\Omega}=0$, which is equivalent to saying that $\mathcal{B}^n_{\Omega_+}\psi -\psi$ vanishes pointwisely in $\partial\Omega$. 
 It is know that $\textbf{H}_\infty$ is self-adjoint and that its spectrum is purely discrete (see \cite{BFSV}, \cite{ATR1} and \cite{OV} for further details). We know also that  $\textbf{H}_M$ tends to $\textbf{H}_\infty$ in strong resolvent sense. In particular, the eigenvalues of $\textbf{H}_M$ converge towards the eigenvalues of $\textbf{H}_\infty$ and any eigenvalue of $\textbf{H}_\infty$ is the limit of eigenvalues of $\textbf{H}_M$, as $M\rightarrow +\infty$ (see \cite{STSV},\cite{RC} and \cite{ATR}). \\
\textbf{Motivations. } In the physics literature, the massive relativistic particles of spin $1/2$ confined in planar or spatial regions are described by Hamiltonian given by the Dirac operator, such systems are of great importance in elementary particle physisc. In planar domain, the operator $\textbf{H}_\infty$ was first considered in 1987 by Berry and Mondragon to study two-dimentional neutrino billard \cite{BM}. Due to its application to model graphene quantum dots (see \cite{AB},\cite{CGPNG}), the operator $\textbf{H}_\infty$ has gained renewed interest from the mathematical physics point of view (\cite{BFSV},\cite{OV},\cite{STSV}). The analogue of this operator in three dimensios was originally introduced by Bogolioubov in the late 60's (see \cite{Bo}) to describe the confinement of quarks in the hadrons. In the mid-70's, this model has been revised by the MIT physicists. With quarks being confined inside a hadron, a useful phenomenological description of quarks in hadrons is provided by the Bag Model (see \cite{CJT}, \cite{CJTW}, \cite{MTY}). While there are many different versions of the model, the MIT Bag Model contains the essential characteristics of the phenomenology of quark confinement. In the MIT bag model, quarks are treated as massless particles inside a bag of finite dimension and infinitely massive outside the bag.
 \\In this present work, we are interested in the asymptotic behavior of the eigenvalues of $\textbf{H}_M$ as $M$ tends to $+\infty$. So our objective is to build  an asymptotic expansion of the eigenvalues of $\textbf{H}_M$ when $M$ tends to $+\infty$. 
 \subsection{Basic notations and main results.} Let us introduce some notations used throughout this article.
 
  \begin{notation} We use the symbols $\langle\cdot ,\cdot\rangle$, $\langle\cdot ,\cdot\rangle_{\Omega_\pm}$ and $(\cdot,\cdot)_{\mathbb{C}^n}$ to denote the scalar products in $L^2(\mathbb{R}^2;\mathbb{C}^2)$ or $L^2(\mathbb{R}^3;\mathbb{C}^4)$, $L^2(\Omega_\pm;\mathbb{C}^n)$ and $\mathbb{C}^n$, respectively. If it is clear from the context we may drop the reference to the spinor space $\mathbb{C}^n$ and simply write $L^2(\mathbb{R}^n)$ and $L^2(\Omega_\pm)$. 
  \end{notation}
  
 \begin{notation}
  Let $\varphi : \Omega_-\longrightarrow \mathbb{R}^+$ be the distance from $x$ to $\partial\Omega$. The open set $\Omega_-$ being smooth, there exits $\Omega^\prime\subset\Omega_-$ a neighborhood of boundary $\partial\Omega$ such that $\varphi$ is smooth on $\Omega^\prime$. Then we have $|\nabla\varphi| =1$ on $\Omega^\prime$ and $\displaystyle\nabla\varphi = \textbf{n}_{\Omega_+}$ on $\partial\Omega$. 
  \end{notation}
  
  \begin{notation} 
  We denote by $ \gamma_0: H^1(\Omega_+)\rightarrow H^{\frac{1}{2}}(\partial\Omega)$ the usual trace operator. Then, for every $\textbf{n}_{\Omega_+}\in \mathbb{S}^{n-1}$ we define the orthogonal projections
  \begin{align}\label{30}
  P_\pm := \frac{1}{2}(1\pm \mathcal{B}^n_{\Omega_+})\gamma_0,
\end{align}  
associated with the eigenvalues $\pm1$ of the matrix $\mathcal{B}^n_{\Omega_+}$.
\end{notation}

\begin{remark} 
 In the 2D case, if one notes by $L$ the length of $\partial\Omega$ and parametrizes $\partial\Omega$ by the curve $\gamma :[0,L]\rightarrow\partial\Omega$ in its arc-length, then \eqref{30} can be written in the following way
\begin{align}\label{31}
P_\pm := \frac{1}{2}(1\pm A(s)), \quad\text{ with }\quad A(s):= \begin{pmatrix}
0& a_2(s)^{*}\\
a_2(s)& 0
\end{pmatrix},\quad a_2(s)^{*}=\bar{a}_2(s)\quad\forall s\in [0,L),
\end{align}
 where $a_2(s):= ie^{i\theta(s)}$, $\theta(s)$ is the angle between $\textbf{n}_{\Omega_+}$ and the $x_1$-axis at the point $\gamma(s)\in\partial\Omega$. 
 In the 3D case the matrix $\mathcal{B}^n_\Omega$ can be written in the form $$\begin{pmatrix}
0 & a_3^\ast\\
a_3 & 0
\end{pmatrix},\quad \text{with }a_3^\ast= -a_3,$$ 
see lemma \ref{appen} for further details, thus for all $z\in\mathbb{R}_+$ we have 
\begin{align}\label{P_+}
P_\pm e^{z\mathcal{B}^n_{\Omega_+}}= e^{\pm z}P_\pm
\end{align}
\end{remark}
 
   \begin{notation}\label{n4}
 For an element $\xi \in \varrho(\textbf{H}_M)\cap\varrho(\textbf{H}_\infty)$, we denote by $R^M(\xi) = (\textbf{H}_M-\xi)^{-1}$ the resolvent of $\textbf{H}_M$, defined on $L^2(\mathbb{R}^n)$ with values in $H^1(\mathbb{R}^n)$ and we write the resolvent of $\textbf{H}_\infty$ by $R^{\infty}(\xi) = (\textbf{H}_\infty-\xi)^{-1}$, defined on $L^2(\Omega)$ with values in $\mathcal{D}^n_\infty$. We also denote by $r_{\pm}:L^2(\mathbb{R}^n)\rightarrow L^2(\Omega_\pm)$ the restriction operator on $\Omega_\pm$ and $e_\pm : L^2(\Omega_\pm)\rightarrow L^2(\mathbb{R}^n)$ the extension by zero in $\Omega_\mp$.  
   \end{notation}
We are now in position to state our main results. We start by showing that the resolvent of $\textbf{H}_M$ admits a first order asymptotic expansion: 

\begin{theorem}\label{th1} The resolvent $R^M(\xi)$ admits an asymptotic expansion of the form:
\begin{align}
R^{M}(\xi) = e_+R^{\infty}(\xi)r_+ + \frac{1}{M}R^{M}_{1}(\xi) + \frac{1}{M}S^{M}_{2}(\xi), \quad \text{ as }\quad M \rightarrow\infty.
\end{align}
Furthermore, there exists a constant $C$ independent of $M$ and $\xi$ such that the operator $R^{M}_{1}(\xi)$ (see the proof of the proposition \ref{pr31} for its definition) and the remainder term $S^M_2$ satifies:
\[\quad\left\{ 
 \begin{array}{l l}
\displaystyle\Vert R^{M}_{1}(\xi)(f)\Vert _{H^2(\Omega_+)}+\frac{1}{\sqrt{M}}\Vert S^{M}_{2}(\xi)(f)\Vert _{H^1(\Omega_+)}\leqslant  C\Vert f\Vert_{H^2(\mathbb{R}^n)}\\ 
\displaystyle \Vert R^{M}_{1}(\xi)(f)\Vert _{L^2(\Omega_-)}+\Vert S^{M}_{2}(\xi)(f)\Vert _{L^2(\Omega_-)}\leqslant  C\Vert f\Vert_{H^2(\mathbb{R}^n)} & \\ \end{array} \right.,\quad\quad\forall f\in H^2(\mathbb{R}^n). \]
\end{theorem}

  Theorem \ref{th1} prove the existence of a first order asymptotic expansion of the resolvent and we deduce the following theorem concerning the eigenvalues :
 
  \begin{theorem}\label{th2}
  Let $\lambda^\infty$ be an eigenvalue of $\textbf{H}_\infty$ with multiplicity $l$. We fix $\eta>0$ such that $B(\lambda^\infty,\eta)\cap \textbf{spec}(\textbf{H}_\infty)=\lbrace\lambda^\infty\rbrace$, where $B(\lambda^\infty,\eta) = \lbrace \xi\in \mathbb{C}, |\lambda^\infty - \xi|\leqslant \eta\rbrace$. Then for sufficiently large $M$, $\textbf{H}_M$ has exactly $l$ eigenvalues $(\lambda^M_k)_{1\leqslant k\leqslant l}$  counted according to their multiplicities in $B(\lambda^\infty,\eta)$ and these eigenvalues admit an asymptotic expansion of the form
  \begin{align}
  \lambda^M_k = \lambda^\infty + \frac{1}{M}\mu_k + \textit{o}(\frac{1}{M}).
  \end{align}
  Where $(\mu_k)_{1\leqslant k\leqslant l}$ are the eigenvalues of the matrix $\mathcal{M}$ with coefficients :
  $$m_{kj}=\frac{1}{2}\langle f_j,f_k\rangle_{\partial\Omega},$$
  with $(f_1,\cdots,f_l)$ an $L^2(\Omega_+)$-orthonormal basis of the eigenspace $\textbf{Ker}(\textbf{H}_\infty -\lambda^\infty I)$.
  \end{theorem}
  
  \begin{remark}
  We can not build an asymptotic expansion of the resolvent at any order in $\mathcal{L}(L^2)$ as it is done at order $1$ in theorem \ref{th1}. Nevertheless it is possible to prove that the eigenvalues of $\textbf{H}_M$ admit an asymptotic expansion at any order.
  \end{remark}
  
  \begin{theorem}\label{th3}
  Let $\lambda^\infty$ be an eigenvalue of $\textbf{H}_\infty$. With the notations of theorem \ref{th2}, for any $k$, $1\leqslant k\leqslant l$, $\lambda^M_k$ admits an asymptotic expansion at any order, that is there exists a sequence $(\mu^j_k)_{j\in\mathbb{N}}$ such that for all $N\in\mathbb{N}$
  \begin{align*}
  \lambda^M_k = \lambda^\infty + \sum_{j=1}^{N}\frac{1}{M^j}\mu_k^j + \textit{O}\left(M^{-N-\frac{1}{2}}\right).
  \end{align*}
  \end{theorem}

  Throughout this paper, we fix an eigenvalue $\lambda^\infty$ of the operator $\textbf{H}_\infty$. Due to the compact embedding of $H^1(\Omega_+)$ in $L^2(\Omega_+)$, we have that the spectrum of $\textbf{H}_\infty$ is discrete, so there exists $\eta>0$ such that $\textbf{spec}(H_\infty)\cap B(\lambda^\infty,\eta)=\lbrace\lambda^\infty\rbrace$.
  \begin{notation}
  We will use the following notations :
  \begin{itemize}
  \item $E^M$ is the sum of the eigenspaces associated with the eigenvalues of $\textbf{H}_M$ contained in $B(\lambda^\infty,\eta)$.
  \item $P^M$ is the spectral projection onto $E^M$. It is given by
  \begin{align}\label{pj1}
  P^M = \frac{-1}{2i\pi}\int_{ \mathcal{C}(\lambda^\infty,\eta)}R^M(\xi) d\xi,
  \end{align}
  where $\mathcal{C}(\lambda^\infty,\eta)$ is the circle of the center $\lambda^\infty$ and radius $\eta$.
  \item $E^\infty$ is the eigenspace with $\textbf{H}_\infty$ associated to $\lambda^\infty$. We denote by $P^\infty$ the spectral projection onto $E^\infty$, which is given by 
  \begin{align}\label{pj2}
  P^\infty = \frac{-1}{2i\pi}\int_{ \mathcal{C}(\lambda^\infty,\eta)}R^{\infty}(\xi) d\xi,
  \end{align}
  \item $\tilde{E}^\infty = e_+(E^\infty)$ and  $\tilde{P}^\infty= e_+P^\infty r_+$, where $e_\pm$ and $r_\pm$ are given by the notation \ref{n4}.
  \end{itemize}
  \end{notation}
  In order to perform an asymptotic expansion of the resolvent we use boundary layers machinery in $\Omega_-$, that is we seek $R^M(\xi)(f)$ of the form:
  $$ R^M(\xi)(f)(x) = V^0(x,M\varphi(x)) + \frac{1}{M}V^1(x,M\varphi(x))+ \frac{1}{M^2}V^2(x,M\varphi(x)) +\cdots,$$
  here $z:=M\varphi(x)$ is the so-called fast variable in penalization methods. We will see during this construction that terms in $e^{-z}$ will appear in the decomposition of the profiles $V^j$ (see equation \eqref{profile} and the subsection~\ref{2.2}), these terms represent the information of the thin layer.
  
 Now we explain briefly the strategy of the proofs of our results. We start by proving a regularity result for the limit operator, we give some estimates and classical results. Next, we show that the resolvent $R^M(\xi)$ admits an asymptotic expansion of order $0$ and we give an estimate for the remainder term with an error of order $1/\sqrt{M}$ as $M\rightarrow\infty$. This proves the convergence of $R^M(\xi)$ to $e_+R^\infty(\xi)r_+$. 
 
 The next step goes as follows: we introduce the free Dirac operator acting in $L^2(\Omega_-)$ with infinite mass boundary condition, in order to avoid the complication of the decomposition of $R^M(\xi)$ in $\Omega_-$. Afterwards, we formally construct the first order asymptotic expansion and we prove Theorem \ref{th1}. In section 3 we prove Theorem \ref{th2}. Following Kato \cite{Kato}, we introduce the operator $\mathcal{U}^M:= I -\tilde{P}^\infty +P^M\tilde{P}^\infty$ which is invertible. Thus the first order asymptotic expansion of $R^M$ gives a first order asymptotic expansion of $\mathcal{U}^M$ in $\mathcal{L}(L^2)$. Afterthat, we introduce $\mathcal{Q}^M= \tilde{P}^\infty(\mathcal{U}^M)^{-1}\textbf{H}_MP^M\mathcal{U}^M\tilde{P}^\infty$ which belongs to $\mathcal{L}(E^\infty)$ and has the same eigenvalues that $\textbf{H}_MP^M$. Using the first order asymptotic expansion of $\mathcal{Q}^M$ and applying classical finite dimensional perturbation theory (see \cite{Kato}) we deduce Theorem \ref{th2}.

 Finally, we denote by $W^M:= (P^M-\tilde{P}^\infty)$ and we introduce the unitary operator $\tilde{\mathcal{U}}^M$ defined by :
 $$\tilde{\mathcal{U}}^M = (I-W^M)^{\frac{-1}{2}})\left( P^M\tilde{P}^\infty + (I-P^M)(I-\tilde{P}^\infty)\right),$$
 we note that $\tilde{\mathcal{Q}}^M= \tilde{P}^\infty(\tilde{\mathcal{U}}^M)^{-1}\textbf{H}_MP^M\tilde{\mathcal{U}}^M\tilde{P}^\infty$ has the same eigenvalues than $\textbf{H}_M$. We will see that it is sufficient to take $f\in C^\infty\cap \mathcal{D}_\infty^n$ and $f=0$ in $\Omega_-$ to build an asymptotic expansion of the resolvent. Since the operator $\tilde{\mathcal{Q}}^M$ is self-adjoint and admit an asymptotic expansion at any order, we conclude that the eigenvalue of $\tilde{\mathcal{Q}}^M$ admit a complete asymptotic expansion. Using the finite dimensional perturabation theory we get Theorem \ref{th3}.

 \subsection{Properties of $\textbf{H}_\infty$}
 In this part, we prove a regularity result and an a priori estimate that will be very useful to estimates the remainder terms.
\begin{theorem}\label{th4}
Let $\Omega_+$ be a domain with $C^{2+k}$-boundary and $f \in H^{k} (\Omega_+ ) $. If $\psi\in\mathcal{D}_{\infty}^n$ is a solution of the following elliptic problem :
\begin{eqnarray}\label{1}
 \left\{
    \begin{split}\
    T_n\psi =  f \quad\text{in } \Omega_+,\\
   P_{-}\psi =  0\quad\text{on } \partial\Omega.   
   \end{split}
  \right.
\end{eqnarray}
   Then $\psi \in H^{k+1}(\Omega_+)$.
\end{theorem}

In order to prove the $H^{k+1}$-regularity, we use a classical resultes to prove the interior regularity and we treat the boundary regularity  separately. For that purpose, we will use tubular coordinates.

\begin{notation}[ Tubular coordinates]\label{note6} Let $X$ be the canonical embedding of $\partial\Omega$ in $\mathbb{R}^n$. We introduce the map $ \kappa : \partial\Omega \times \left(-\delta, \delta \right)\rightarrow \mathbb{R}^n$ defined by 
$$ \kappa (s,t) = X(s) +t\textbf{n}_{\Omega_+}(s).$$
This transformation is a $C^1$-diffeomorphism for $\delta\in (0, \delta_0)$ provided that $\delta_0$ is sufficiently small. For $\delta<\delta_0$, we define a neighbourhood $Q_\delta$ of $\partial\Omega$ by
$$ Q_\delta : = \lbrace x\in \mathbb{R}^n : \text{dist} (x,\partial\Omega ) < \delta \rbrace.$$
Then $\kappa$ is a local parametrization of $Q_\delta$ for $\delta_0$ small enough.
\end{notation}

\textbf{Proof of Theorem \ref{th4} . }Let us first show this result in the $2$D case.
By induction, assume that $k=1$ i.e $\partial\Omega \in C^3$ and $f\in H^{1}(\Omega_+)$. Applying $T_2$ to \eqref{1},  we get that $\displaystyle -\Delta\psi =T_2 f$ holds, in a distributional sense. Using \cite[Theorem 8.8]{Gil}, one gets that for any $\delta<\delta_0$ we have that $\psi\in H^2(\Omega_\delta)$, where $\Omega_\delta := \lbrace x\in\Omega_+ : \text{dist}(x, \partial\Omega )> \delta/2 \rbrace$. It remains to show that $\psi\in H^2(U_\delta)$, where $U_\delta:=Q_\delta\cap \Omega_+$. For this we define the following operator in $L^2(Q_\delta;\mathbb{C}^2)$
\begin{align}
(P\psi)(x) =\frac{1}{2}(1-A(s))\psi(\kappa(s,t)),\quad x=\kappa(s,t)\in U_\delta,
\end{align}
where $A$ is the matrix function defined in \eqref{31}. Let $x_0 \in\partial\Omega$ and $\chi\in C^\infty(\mathbb{R}^2,[0,1])$ supported in $B(x_0,2\rho)$ with $\chi=0$ on $\mathbb{R}^2 \setminus B(x_0,2\rho)$ and $\chi =1$ in $B(x_0,\rho)$. We choose $\displaystyle \rho<\delta/2$. Analysis similar to that in the proof of \cite[Theorem 2]{STSV} shows
that 
\begin{eqnarray}\label{2}
 \left\{
    \begin{split}\
    \displaystyle -\Delta(P\chi \psi) =  g \quad\text{in } \Omega_+,\\
   P\chi\psi =  0\quad\text{on } \partial\Omega,   
   \end{split}
  \right.
\end{eqnarray}
holds for $g= T_2(P\chi f+\frac{1}{2}[T_2,\chi]\psi) -\frac{1}{2}T_2[T_2,A\chi]\psi \in L^2(\Omega_+)$. As a consequence we get that $P\chi \psi\in H^{2}(\Omega_+)$ and $P\chi \psi =0$ on $\partial\Omega$ (see \cite[Theorem 8.11]{Gil}). Combine this with the  compactness of the boundary, it provides that 
\begin{align}
\psi_2 -ie^{i\theta}\psi_1 \in H^2(U).
\end{align}
It follows from this that
\begin{align}
-i(\partial_1-i\partial_2)\psi_2 - e^{i\theta}(\partial_1 -i\partial_2)\psi_1 \in H^1(U_\delta).
\end{align}
Since $-i(\partial_1-i\partial_2)\psi_2=f_1\in H^1(U_\delta)$, we find that 
\begin{align}
(\partial_1 -i\partial_2)\psi_1\in H^1(U_\delta).
\end{align}
As a consequence, we get  
\begin{align*}
(\partial_1 +i\partial_2)(\partial_1 -i\partial_2)\psi_1= \Delta\psi_1\in L^2(U_\delta),\\
(\partial_1 -i\partial_2)(\partial_1 -i\partial_2)\psi_1= (\partial_1^{2}-\partial_2^{2}-2i\partial_{12})\psi_1	\in L^2(U_\delta),\\
\partial_1(\partial_1-i\partial_2)\psi_1=(\partial_{1}^{2}-i\partial_{12})\psi_1 \in L^2(U_\delta).
\end{align*}
Combine the last equalities, we get that $\partial_{2}^{2}\psi_1$ and  $\partial_{12}\psi_1$ belongs to $L^2(U_\delta)$. Since $\Delta \psi_1\in L^2(U_\delta)$ we obtain also that $\partial_1^{2}\psi_1\in L^2(U_\delta)$. That $\psi_1\in H^2(\Omega_+)$ follows from this and interior regularity. The same conclusion can be drawn for $\psi_2$, i.e $\psi_2\in H^2(\Omega_+)$. This completes the proof for $k=1$.
Assume now that  $\partial\Omega \in C^{2+k}$, $f\in H^{k}(\Omega_+)$ and $\psi\in\mathcal{D}_{\infty}^2\cap H^{k}(\Omega_+)$. Using  \cite[Theorem 8.10]{Gil}, we get that $\psi \in H^{k+1}(\Omega_\delta)$ for some $\delta<\delta_0$.  Furthermore we get that \eqref{2} holds for $g= T_2(P\chi f+\frac{1}{2}[T_2,\chi]\psi) -\frac{1}{2}T_2[T_2,A\chi]\psi \in H^{k-1}(\Omega_+)$ and equation \eqref{2} implies by \cite[Theorem 8.13]{Gil} that $P\chi\psi\in H^{k}(\Omega_+)$ and $P\chi\psi = 0$ on $\partial\Omega$. As a consequence we get
\begin{align*}
\psi_2 -ie^{i\theta}\psi_1\in H^{k}(U_\delta).
\end{align*}
The rest of the proof runs as before. For the 3D case we give only the main ideas of the proof. Let $k=1$. Denote by $\theta$ the azimuthal angle in the $xy$-plane from the $x$-axis with $0\leqslant\theta<2\pi$ and $\phi$ the polar angle from the positive $z$-axis with $0\leqslant\phi\leqslant\pi$.

 We write the outward normal vector in the spherical coordinates, i.e :
  $$\textbf{n}_{\Omega_+}(s) =(\sin\phi(s)\cos\theta(s),\sin\phi(s)\sin\theta(s),\cos\phi(s)).$$   Then the same method as in the $2$D case  provides that $P\chi\varphi \in H^{2}(\Omega_+)$ and $P\chi\varphi =0$ on $\partial\Omega$. As consequence we obtain that 
\begin{align}
&\label{I}\psi_1 + i\cos\phi\psi_3 +i\sin\phi e^{-i\theta}\psi_4 \in H^2(U_\delta),\\
& \label{II}\psi_2 + i\sin\phi e^{i\theta}\psi_3 -i\cos\phi \psi_4 \in H^2(U_\delta),\\
& \label{III}i\cos\phi \psi_1 +i\sin\phi e^{-i\theta}\psi_2 +\psi_3 \in H^2(U_\delta),\\
& i\sin\phi e^{i\theta}\psi_1 -i\cos\phi \psi_2 + \psi_4 \in H^2(U_\delta). 
\end{align}
Since $\partial\Omega \in C^3$ we have that $\theta, \phi \in C^2$. Then $\eqref{III} -i\cos\phi \eqref{I} - i\sin\phi e^{-i\theta}\eqref{II} = 2\psi_3 \in H^2(U_\delta)$. With the same reasoning we show that $\psi_j \in H^2(U_\delta)$, $j=1,..,4$. We conclude with the same arguments as in the 2D case. \qed \\
In the remainder of this work we frequently use the following Trace Theorem (see \cite{A}).
\begin{theorem} There exists $C>0$ such that
\begin{align*}
\Vert \psi\Vert_{L^2(\partial\Omega)} \leqslant C \Vert \psi\Vert_{H^1(\Omega_+)},\quad \forall \psi\in H^1(\Omega_+).
\end{align*}
\end{theorem}
By duality we easily obtain the following lemma which will be very useful afterwards to estimate the remainder terms.
\begin{lemma}\label{lemme 1} There exists a constant $C$ such that for any $\xi \in\mathcal{C}(\lambda^\infty,\eta)$, for any $\psi\in H^1(\Omega_+)$ such that
$(T_n-\xi)\psi\in L^2(\Omega_+)$, we have
\begin{align*}
\Vert \psi\Vert_{L^2(\Omega_+)} \leqslant C\big(\Vert (T_n-\xi)\psi\Vert_{L^2(\Omega_+)} +\Vert P_- \psi\Vert_{L^2(\partial\Omega_+)}\big).
\end{align*}
\end{lemma}
\textbf{Proof. }Let $f\in L^2(\mathbb{R}^n)$. There exists $g\in \mathcal{D}^n_\infty$ such that $(T_n-\xi)g=f$ in $\Omega_+$. Furthermore $P_-g=0$ and $||g||_{H^1(\Omega_+)}\leqslant C||f||_{L^2(\Omega_+)}$, where $C$ does not depend on $f$ and $\xi$. Then,  we have :
\begin{align*}
 \int_{\Omega_+} (\psi, f)_{\mathbb{C}^n}  &= \int_{\Omega_+}(\psi, (T_n-\xi)g)_{\mathbb{C}^n} \\
&=\int_{\Omega_+}(T_n\psi,g)_{\mathbb{C}^n}  -\xi\int_{\Omega_+}(\psi, g)_{\mathbb{C}^n}- \langle\psi , \Theta_n \mathcal{B}^n_{\Omega_+} g\rangle _{\partial\Omega}\\
&=\int_{\Omega_+}((T_n-\xi)\psi, g)_{\mathbb{C}^n} - \langle\psi , \Theta_n \mathcal{B}^n_{\Omega_+} P_+g\rangle _{\partial\Omega}\\
&=\int_{\Omega_+}((T_n-\xi)\psi ,g)_{\mathbb{C}^n} - \langle\Theta_n P_-\psi ,  g\rangle _{\partial\Omega}.
\end{align*}
Here we used some properties of the matrices $\mathcal{B}^n_{\Omega_+}$ and $P_\pm$ (see proposition \ref{pro ape}) and Green's Formula. Thus
\begin{align*}
\left|\int_{\Omega_+} (\psi, f)_{\mathbb{C}^n}\right| &\leqslant \Vert (T_n-\xi)\psi\Vert_{L^2(\Omega_+)}\Vert g\Vert_{L^2(\Omega_+)} + \Vert P_- \psi\Vert_{L^2(\partial\Omega)}\Vert g\Vert_{L^2(\partial\Omega)}\\
&\leqslant C\big(\Vert (T_n-\xi)\psi\Vert_{L^2(\Omega_+)} +\Vert P_- \psi\Vert_{L^2(\partial\Omega)}\big)\Vert g\Vert_{H^1(\Omega_+)}\\
 &\leqslant C\big(\Vert (T_n-\xi)\psi\Vert_{L^2(\Omega_+)} +\Vert P_- \psi\Vert_{L^2(\partial\Omega)}\big)\Vert f\Vert_{L^2(\Omega_+)}.
\end{align*}
where in the last equality, we use the trace theorem for $g$ and the inequality $||g||_{H^1(\Omega_+)}\leqslant C||f||_{L^2(\Omega_+)}$. Since the estimates are true for all $f\in L^2(\Omega_+)$,  we obtain that
\begin{align*}
\Vert \psi\Vert_{L^2(\Omega_+)} \leqslant C\big(\Vert (T_n-\xi)\psi\Vert_{L^2(\Omega_+)} +\Vert P_- \psi\Vert_{L^2(\partial\Omega)}\big).
\end{align*}\qed

\section{Asymptotic expansion of the resolvent}
We follow the same strategy used for Schr\"{o}dinger operator in \cite{BrCa02}.\label{sect} So we fix $f\in L^2(\mathbb{R}^n)$ and we denote by $u^M$ (resp $v^M$) the restriction of $R^M(\xi)(f)$ on $\Omega_+$ (resp $ \Omega_-$). Thus the couple $(u^M, v^M)$ satisfies the following system 
\begin{equation}\label{23}
  \left\{
    \begin{aligned}
      (T_n-\xi) u^M&=f&\text{ in }\text{ }\Omega_+ \\
       (T_n-\xi) v^M + M\Theta_n v^M &=f&\text{ in }\text{ }\Omega_- \\
       u^M &= v^M& \text{ on }\text{ } \partial\Omega \\ 
    \end{aligned}
  \right.
\end{equation}
We will seek an asymptotic expansion of $u^M$ and $v^M$ of the form :
\vspace{-1em}
\begin{align*}
&u^M(x) = U^0(x) +\frac{1}{M}U^{1 }(x) +\frac{1}{M^2}U^{2}(x)+\ldots\quad\quad\quad\quad\quad\quad\quad\quad\quad\quad\quad\quad\quad\text{ }\\
&v^M(x)= V^0(x,M\varphi(x))+ \frac{1}{M}V^{ 1}(x,M\varphi(x))+\frac{1}{M^2}V^{ 2}(x,M\varphi(x))+\ldots,
\end{align*}
where the profiles $V^j(x,z)$ can be decomposed on the following way:
\begin{align}\label{profile}
V^j(x,z) = \tilde{V}^j(x,z) +\overline{V}^j(x),
\end{align}
here $\overline{V}^j(x)=\lim\limits_{z\longrightarrow\infty}V^j(x,z)$ and we assume that $\tilde{V}^j$, as well as all its partial derivatives in both variables $x$ and $z$ tend to zero when $z$ tends to $+\infty$.
In order to avoid the complication of this decomposition, we introduce for $M>0$ the self-adjoint operator acting in $L^2(\Omega_-)$ defined by:
\begin{align}
\textbf{H}^{M}_{+}\psi = T_n\psi +M\Theta_n\psi,\quad\quad  \psi\in \mathcal{D}_+^n :=\lbrace \psi\in H^1( \Omega_-) :P_+\psi =0 \text{ on } \partial\Omega\rbrace.
\end{align} 
For $f\in L^2(\mathbb{R}^n)$ we define $\omega^M(x) := (\textbf{H}^{M}_{+}-\xi)^{-1}r_-f$, which admits an asymptotic expansion (see the proof of theorem \ref{th1}) on the form :
 \begin{align*}
 \omega^M = \omega^{0}+\frac{1}{M}\omega^{1}+\frac{1}{M^2}\omega^{2}+\ldots., \text{ in } L^2(\Omega_-).
 \end{align*}
Note that when $x$ is far from the boundary $\partial\Omega$, the equality $\overline{V}^j= \omega^j$ holds. We define the operator $\tilde{R}^M(\xi)(f) : = \left( (\textbf{H}_M-\xi)^{-1}-0\oplus(\textbf{H}^{M}_{+}-\xi)^{-1}\right) $, thus $\tilde{u}^M:=u^M$ is the restriction of $\tilde{R}^M(\xi)(f)$ on $\Omega_+$ and $\tilde{v}^M:=v^M -\omega^M$ is the restriction of $\tilde{R}^M(\xi)(f)$ on $\Omega_-$. In this way we obtain that the couple $(\tilde{u}^M,\tilde{v}^M)$ satisfies the following system :
\begin{equation}\label{3} 
  \left\{
    \begin{aligned}
      (T_n-\xi) \tilde{u}^M&=f&\text{ in }\text{ }\Omega_+\quad\quad (1)\\
       (T_n-\xi) \tilde{v}^M + M\Theta_n \tilde{v}^M &=0&\text{in }\text{ }\Omega_- \quad\quad (2)\\
       P_+\tilde{u}^M &= P_+\tilde{v}^M& \text{ on }\text{ } \partial\Omega\quad\quad (3) \\ 
        P_-\tilde{u}^M &= P_-\tilde{v}^M +P_-\omega^M& \text{ on }\text{ } \partial\Omega\quad\quad (4) \\
    \end{aligned}
  \right.
\end{equation}   
We replace formally $\tilde{u}^M$ and $\tilde{v}^M$ by their profile and we identify the different terms in the power of $M$.
\subsection{Order 0}
We begin by recalling a statement that follows from \cite[Lemma 2]{STSV}, \cite[Lemma 4]{STSV} and \cite[Lemma 1.6]{ATR}.
Let $\psi \in H^1(\mathbb{R}^n)$, a  direct computation shows that
\begin{align}
 \displaystyle \Vert (T_n+M\Theta_n)\phi\Vert ^{2}_{L^2(\Omega_-)} &= \Vert \nabla\phi \Vert^2_{L^2(\Omega_-)}+ M^2\Vert \phi \Vert^2_{L^2(\Omega_-)}-M\Vert P_+\phi\Vert^{2}_{L^2(\partial\Omega)}+ M\Vert P_-\phi\Vert^{2}_{L^2(\partial\Omega)}\label{41}
\end{align}
 The estimation of the first three terms above, gives
\begin{proposition}[\cite{STSV},\cite{ATR}]  There exist constants $c$ and $M_0>1$ such that for all $\phi\in H^1(\mathbb{R}^n)$ and $M>M_0$ holds 
\begin{align}
  \displaystyle \Vert (T_n+M\Theta_n)\phi\Vert ^{2}_{L^2(\Omega_-)}& \geqslant  M\Vert P_-\phi\Vert^{2}_{L^2(\partial\Omega)}-c\Vert \phi\Vert^{2}_{L^2(\partial\Omega)}.\label{4}
\end{align}
\end{proposition}
\begin{proposition}\label{prop11} Let $\psi \in H^1(\Omega_+)$, $\phi\in H^1(\Omega_-)$ and $f\in L^2(\mathbb{R }^n)$ such that 
\begin{equation}
  \left\{
    \begin{aligned}
      \displaystyle  (T_n-\xi)\psi& = 0\quad\quad\quad\text{ in }\text{  }\Omega_+  \\
       \displaystyle  (T_n-\xi)\phi+ M\Theta_n\phi &= f\quad\quad\quad\text{ in }\text{  }\Omega_- \\
     \displaystyle P_+(\psi+\gamma) &= P_+\phi\quad\text{ on }\text{  } \partial\Omega \\ 
    \end{aligned}
  \right.
\end{equation}
for $\gamma\in H^1(\Omega_+)$. Assume that
\begin{align} 
\Vert P_-\phi\Vert^{2}_{L^2(\partial\Omega)} &\geqslant \Vert P_-\psi\Vert^{2}_{L^2(\partial\Omega)} - \frac{C}{M}\Vert f\Vert^{2}_{L^2(\mathbb{R}^n)}\label{P}.
\end{align} 
         Then there exist $C>0$ such that for $M$ large enough, we have
\begin{align*}
\Vert \phi\Vert^{2}_{L^2(\Omega_-)}&\leqslant \frac{C}{M}\left(\Vert f\Vert^{2}_{L^2(\mathbb{R}^n)}+M\Vert P_+\gamma\Vert^{2}_{L^2(\partial\Omega)}\right),\\
\Vert \phi\Vert^{2}_{H^1(\Omega_-)}&\leqslant C\left(\Vert f\Vert^{2}_{L^2(\mathbb{R}^n)}+M\Vert P_+\gamma\Vert^{2}_{L^2(\partial\Omega)}\right),\\
\Vert \psi\Vert^{2}_{H^1(\Omega_+)}&\leqslant \frac{C}{M}\left(\Vert f\Vert^{2}_{L^2(\mathbb{R}^n)}+\Vert P_+\gamma\Vert^{2}_{L^2(\partial\Omega)}\right).
\end{align*}
In particular, if $P_+\phi=P_+\psi$, then 
\begin{align*}
M\Vert \psi\Vert^{2}_{H^1(\Omega_+)}+M^2\Vert \phi\Vert^{2}_{L^2(\Omega_-)}&\leqslant C\Vert f\Vert^{2}_{L^2(\mathbb{R}^n)},\\
\Vert \phi\Vert^{2}_{H^1(\Omega_-)}&\leqslant C\Vert f\Vert^{2}_{L^2(\mathbb{R}^n)}.
\end{align*}
\end{proposition}
\textbf{Proof.} Using the inequality \eqref{4}  and the equation satisfied by $\phi$, we get
\begin{align*}
 \displaystyle M\Vert P_-\phi\Vert^{2}_{L^2(\partial\Omega)}\leqslant C\left( |\xi|^2 \Vert \phi\Vert ^{2}_{L^2(\Omega_-)}+\Vert f \Vert^2_{L^2(\Omega_-)}+\Vert P_+\phi\Vert^{2}_{L^2(\partial\Omega)} \right).
\end{align*}
Thus, using the boundary condition and the inequality \eqref{P}, we obtain
\begin{align*}
 \displaystyle M\Vert P_-\psi\Vert^{2}_{L^2(\partial\Omega)}\leqslant C\left( |\xi|^2 \Vert \phi\Vert ^{2}_{L^2(\Omega_-)}+\Vert f \Vert^2_{L^2(\Omega_-)}+\Vert P_+(\psi+\gamma)\Vert^{2}_{L^2(\partial\Omega)} \right).
\end{align*}
Then, apply lemma \ref{lemme 1} to $\psi$, it provides 
\begin{align*}
 \displaystyle \Vert \psi\Vert^{2}_{L^2(\Omega_+)}\leqslant \frac{C}{M}\left( |\xi|^2 \Vert \phi\Vert ^{2}_{L^2(\Omega_-)}+\Vert f \Vert^2_{L^2(\Omega_-)}+\Vert P_+\psi\Vert^{2}_{L^2(\partial\Omega)}+\Vert P_+\gamma\Vert^{2}_{L^2(\partial\Omega)} \right).
\end{align*}
Since $ \Vert T_n \psi\Vert_{L^2(\Omega_+)}=\Vert \nabla \psi\Vert_{L^2(\Omega_+)}=|\xi|\Vert \psi\Vert_{L^2(\Omega_+)} $, so we get 
\begin{align*}
 \displaystyle \Vert \psi\Vert^{2}_{H^1(\Omega_+)}\leqslant \frac{C(1+|\xi|^2)}{M}\left( |\xi|^2 \Vert \phi\Vert ^{2}_{L^2(\Omega_-)}+\Vert f \Vert^2_{L^2(\Omega_-)}+\Vert P_+\psi\Vert^{2}_{L^2(\partial\Omega)}+\Vert P_+\gamma\Vert^{2}_{L^2(\partial\Omega)} \right).
\end{align*}
 Hence for $M$ large enough, applying the trace theorem to $P_+\psi$, one gets 
\begin{align}\label{psi}
\Vert \psi\Vert^{2}_{H^1(\Omega_+)}&\leqslant \frac{C}{M}\left(\Vert f\Vert^{2}_{L^2(\mathbb{R}^n)}+ |\xi|^2\Vert \phi\Vert^{2}_{L^2(\Omega_-)}+\Vert P_+\gamma\Vert^{2}_{L^2(\partial\Omega)} \right).
\end{align}
Remark that $(T_n+M\Theta_n)\phi = \xi\phi +f$, then from the equality  \eqref{41} and the boundary condition we get 
\begin{align*}
 \displaystyle \Vert \nabla\phi \Vert^2_{L^2(\Omega_-)}+ M^2\Vert \phi \Vert^2_{L^2(\Omega_-)} \leqslant \Vert f\Vert ^{2}_{L^2(\Omega_-)} + |\xi|^2\Vert \phi\Vert^{2}_{L^2(\Omega_-)}+M\Vert P_+(\psi+\gamma)\Vert^{2}_{L^2(\partial\Omega)}
\end{align*}
Therefore using the trace theorem in the last inequality and injecting inequality \eqref{psi} into, it provides
\begin{align}
 \displaystyle \Vert \nabla\phi \Vert^2_{L^2(\Omega_-)}+ M^2\Vert \phi \Vert^2_{L^2(\Omega_-)} \leqslant C\left(\Vert f\Vert ^{2}_{L^2(\Omega_-)} + |\xi|^2\Vert \phi\Vert^{2}_{L^2(\Omega_-)}+M\Vert P_+\gamma\Vert^{2}_{L^2(\partial\Omega)}\right).\label{2.9}
\end{align}
From this we obtain 
\begin{align*}
\Vert \phi\Vert^{2}_{L^2(\Omega_-)}&\leqslant \frac{C}{M^2}\left(\Vert f\Vert^{2}_{L^2(\mathbb{R}^n)}+M\Vert P_+\gamma\Vert^{2}_{L^2(\partial\Omega)}\right),\\
\Vert \phi\Vert^{2}_{H^1(\Omega_-)}&\leqslant C\left(\Vert f\Vert^{2}_{L^2(\mathbb{R}^n)}+M\Vert P_+\gamma\Vert^{2}_{L^2(\partial\Omega)}\right).
\end{align*}
Using this and inequality \eqref{psi}, we get
$$\Vert \psi\Vert^{2}_{H^1(\Omega_+)}\leqslant \frac{C}{M}\left(\Vert f\Vert^{2}_{L^2(\mathbb{R}^n)}+ \Vert P_+\gamma\Vert^{2}_{L^2(\partial\Omega)}\right),$$
In particular if $P_+\gamma = 0$, inequality \eqref{2.9} provides 
\begin{align*}
M\Vert \psi\Vert^{2}_{H^1(\Omega_+)}+M^2\Vert \phi\Vert^{2}_{L^2(\Omega_-)}&\leqslant C\Vert f\Vert^{2}_{L^2(\mathbb{R}^n)},\\
\Vert \phi\Vert^{2}_{H^1(\Omega_-)}&\leqslant C\Vert f\Vert^{2}_{L^2(\mathbb{R}^n)}.
\end{align*}
This gives the desired result. \qed
\begin{remark}\label{rem2} Note that if there exists $C>0$ such that: $\Vert P_+\gamma\Vert_{L^2(\partial\Omega)}\leqslant C\Vert f\Vert^{2}_{L^2(\mathbb{R}^n)}$, then 
\begin{align*}
\Vert \phi\Vert^{2}_{L^2(\Omega_-)}+\Vert \psi\Vert^{2}_{H^1(\Omega_+)}\leqslant \frac{C}{M}\Vert f\Vert^{2}_{L^2(\mathbb{R}^n)}.
\end{align*}
\end{remark}
\begin{proposition}
The resolvent $R^M(\xi)$ admits the following asymptotic expansion:
\begin{align*}
 R^M(\xi) = e_+R^\infty(\xi)r_+ + S^{M}_{1}(\xi),
\end{align*}
 where  $S^{M}_{1}(\xi)$ is an  operator on $L^2(\mathbb{R}^n)$. Moreover, there exists a constant $C$ such that 
 \begin{align*}
   \displaystyle\Vert S^{M}_{1}(\xi)(f)\Vert _{H^1(\Omega_+)}+\Vert S^{M}_{1}(\xi)(f)\Vert _{L^2(\Omega_-)}\leqslant  \frac{C}{\sqrt{M}}\Vert f\Vert_{L^2(\mathbb{R}^2)},\quad\forall M>M_0,\forall f\in L^2(\mathbb{R}^n)
  \end{align*}
\end{proposition}
 
\textbf{Proof. } Let $f\in L^2(\mathbb{R}^n)$, we know that $u^M$ and $v^M$ satisfy \eqref{23}. We write
 \[  \left\{ 
\begin{array}{l l}
  \displaystyle u^M(x) = U^0(x) + a^M(x),\\
  \displaystyle v^M(x)= b^M(x), & \\ \end{array} \right. \]

where $U^0$ satisfies : 
\[  \left\{ 
\begin{array}{l l}
  \displaystyle  (T_n-\xi)U^0 = f\quad\text{in }\Omega_+  ,&\\
  \displaystyle P_-U^0 = 0\quad\text{on } \partial\Omega, & \\ \end{array} \right. \]
thus 
$$ U^0= e_+R^\infty(\xi)r_+(f).$$
The remainder terms satisfy the following system :
\begin{equation*}
\left\{
    \begin{aligned}
      (T_n-\xi)a^M& = 0\text{ in }\text{ }\Omega_+ \\
       (T_n-\xi)b^M+ M\Theta_nb^M &= f\text{ in }\text{ }\Omega_- \\
       P_+a^M +P_+U^0 &= P_+b^M \text{ on }\text{ } \partial\Omega \\ 
       P_-a^M  &= P_-b^M \text{ on }\text{ } \partial\Omega \\
    \end{aligned}
  \right.
  \end{equation*}
  It is noted that there exist a constants $C_1$ and $C_2$  independent on $\xi$ and $f$ such that:
\begin{align}\label{5}
\Vert TU^0\Vert^{2}_{L^2(\Omega_+)} + \Vert U^0\Vert^{2}_{L^2(\Omega_+)}&\leqslant C_1 \Vert f\Vert^{2}_{L^2(\Omega_+)},\\
\label{6}\Vert U^0\Vert^{2}_{L^2(\partial\Omega)}&\leqslant C_2\left( \Vert TU^0\Vert^{2}_{L^2(\Omega_+)} +\Vert U^0\Vert^{2}_{L^2(\Omega_+)}\right),
\end{align}
here we used the trace theorem. In addition, boundary conditions give us  
\begin{align*}
\Vert P_-b^M\Vert^{2}_{L^2(\partial\Omega)} \geqslant \Vert P_-a^M\Vert^{2}_{L^2(\partial\Omega)}- \frac{C}{M}\Vert f\Vert^{2}_{L^2(\mathbb{R}^n)}.
\end{align*} 
Furthermore, inequalities \eqref{5} and \eqref{6} implies
\begin{align*}
 \Vert P_+U^0\Vert^{2}_{L^2(\partial\Omega_+)}&\leqslant C_1 \Vert f\Vert^{2}_{L^2(\Omega_+)}.
\end{align*}
Hence, apply proposition \ref{prop11} with $\gamma = U^0$ and taking remarque \ref{rem2} into account, one obtains
\begin{align*}
\Vert a^M\Vert^{2}_{H^1(\Omega_+)}&\leqslant \frac{C}{M}\Vert f\Vert^{2}_{L^2(\mathbb{R}^n)},\\
\Vert b^M\Vert^{2}_{L^2(\Omega_-)}&\leqslant \frac{C}{M}\Vert f\Vert^{2}_{L^2(\mathbb{R}^n)}.
\end{align*}
Now we have 
\[ S^{M}_{1}(\xi)= \left\{ 
\begin{array}{l l}
  \displaystyle a^{M}(x)\text{ if } x\in \Omega_+,\\
  \displaystyle b^{M}(x)\text{ if } x\in\Omega_-. & \\ \end{array} \right. \]
This completes the proof. \qed
\subsection{First order asymptotic expansion}
Using the previous notations, we will seek an asymptotic expansion of $\tilde{u}^M$ and $\tilde{v}^M$ of the form:\label{2.2}
\begin{align*}
\tilde{u}^M(x) &= U^0(x) +\frac{1}{M}U^{1 }(x) +\frac{1}{M^2}U^{2}(x)+\ldots\quad\quad\quad\quad\quad\quad\quad\quad\quad\quad\quad\quad\quad\text{ }\\
\tilde{v}^M(x)&= V^0(x,M\varphi(x))+ \frac{1}{M}V^{ 1}(x,M\varphi(x))+\frac{1}{M^2}V^{ 2}(x,M\varphi(x))+\ldots
\end{align*}
Where $(\tilde{u}^M,\tilde{v}^M)$ satisfy \eqref{3} and $V^i=\sum_{j\geqslant1}(M\varphi)^{j-1}\tilde{\alpha}_{i,j}e^{-M\varphi}$, for some functions $\tilde{\alpha}_{i,j}$ defined from the boundary. We will replace $\tilde{u}^M$ and $\tilde{v}^M$ by their asymptotic expansion in \eqref{3} and we will identify the different powers of $M$. We remark that 
\begin{align*}
T_2 V^i(x,M\varphi(x))&=-\textit{i}\begin{pmatrix}
0 & \partial_1 -\textit{i} \partial_2\\
\partial_1 +\textit{i}\partial_2 & 0
\end{pmatrix}
\begin{pmatrix}
V^{i}_{1}(x,M\varphi(x))\\
V^{i}_{2}(x,M\varphi(x))
\end{pmatrix}\\
&= \begin{pmatrix}
-\textit{i}\partial_1V^{i}_{2}-\textit{i}M\partial_1\varphi V^{i}_{z,2} -M\partial_2\varphi V^{i}_{z,2}-\partial_2V^{i}_{2}\\
-\textit{i}\partial_1 V^{i}_{1}-\textit{i}M\partial_1\varphi V^{i}_{z,1}+M\partial_2\varphi V^{i}_{z,1}+\partial_2V^{i}_{1}
\end{pmatrix}\\
&=T_{2,x}V(x,M\varphi(x)) + M (T_2\varphi(x))V^{i}_{z}(x,M\varphi(x)),
\end{align*}
and we get  the same thing with $T_3$, where we use for $V^{j}=(V^{j}_{1},\ldots,V^{j}_{n})^t
$ the notations $\displaystyle V_z =\frac{\partial V}{\partial z}$ and $T_{n,x}$ acting as a derivative on the $x$ variable. Using that $|\nabla\varphi |=1$ near the boundary, we formally obtain the following equations :

\[ \text{from (2.3.1)}, \textbf{ }\left\{ 
\begin{array}{l l}
  (\rm{i})\quad\quad\displaystyle &(T_n-\xi)U^0=f\quad\text{ in } \Omega_+,\\
  (\rm{ii})\quad\quad\displaystyle &(T_n-\xi)U^1=0\quad\text{ in } \Omega_+,  \\ \end{array} \right. \]

\[ \text{from (2.3.2)}, \textbf{ }\left\{ 
\begin{array}{l l}
  (\rm{iii})\quad\quad\displaystyle &V^{0}_{z}=-\mathcal{B}^n_{\Omega_+} V^0\quad\text{ in } \Omega_-,\\
  (\rm{iv})\quad\quad\displaystyle &(T_{n,x}-\xi)V^0 +\Theta_n\mathcal{B}^n_{\Omega_+} V^{1}_{z}+\Theta_nV^1=0\quad\text{ in }\Omega_-,\\
   (\rm{v})\quad\quad\displaystyle &(T_{n,x}-\xi)V^1 +\Theta_n\mathcal{B}^n_{\Omega_+} V^{2}_{z}+\Theta_nV^2=0\quad\text{ in } \Omega_-, \\ \end{array} \right. \]

\[ \text{from (2.3.3)}, \textbf{ }\left\{ 
\begin{array}{l l}
  (\rm{vi})\quad\quad\displaystyle &P_+U^0=P_+V^0\quad\text{ on } \partial\Omega,\\
  (\rm{vii})\quad\quad\displaystyle &P_-U^0=P_-V^0\quad\text{ on } \partial\Omega,\\
  (\rm{viii})\quad\quad\displaystyle &P_+U^1=P_+V^1\quad\text{ on } \partial\Omega,  \\ \end{array} \right. \]
here we used the property \eqref{P_+} and the fact that $ \Theta_n\mathcal{B}^n_{\Omega_+}= T_n\varphi$ (it's just a simple computation).
  Using the last equations,  we obtain by the Duhamel formula:
\begin{itemize}
\item{ from (\rm{iii}), we have $\displaystyle V^0(x,z)= e^{-z\mathcal{B}^n_{\Omega_+}}\alpha_{0}(x)$, where $\alpha_0=V^{0}_{|z=0}$.}
\item{from (\rm{iv}), $\displaystyle V^{1}(x,z) = e^{-z\mathcal{B}^n_{\Omega_+}}\alpha_{1}(x) - \int^{z}_{0}e^{-(z-s)\mathcal{B}^n_\Omega}\mathcal{B}^n_{\Omega_+}\Theta_n (T_{n,x}-\xi)V^0(x,s)\text{ds}$, where $\alpha_1=V^{1}_{|z=0}$.}
\item{from (\rm{iv}), $\displaystyle V^{2}(x,z) = e^{-z\mathcal{B}^n_{\Omega_+}}\alpha_{2}(x) - \int^{z}_{0}e^{-(z-s)\mathcal{B}^n_{\Omega_+}}\mathcal{B}^n_\Omega\Theta_n (T_{n,x}-\xi)V^1(x,s)\text{ds}$, where $\alpha_2=V^{2}_{|z=0}$.}
\end{itemize}
Since $\tilde{u}^M$ and $\tilde{v}^M$ are bounded, this imposes the following conditions:
\[ \textbf{ }\left\{ 
\begin{array}{l l}
  (\rm{ix})\quad\quad\displaystyle &P_-V^0=0\quad\text{ on } \partial\Omega,\\
  (\rm{x})\quad\quad\displaystyle& P_-V^1=\emph{O}(1),\\
  (\rm{xi})\quad\quad\displaystyle &P_-V^2=\emph{O}(1). \\ \end{array} \right. \]
  \begin{itemize}
\item{ Since (\rm{ix}) implies that $P_-U^0 = 0$ on $\partial\Omega$, using (\rm{i}) and (\rm{vii}) we get:}
  \begin{equation}\label{10} 
  \left\{
  \begin{aligned}
 \displaystyle (T_n-\xi)U^0&=f\quad\quad\text{ in } \Omega_+,\\
   \displaystyle P_-U^0&= 0 \quad\text{ on } \partial\Omega.
  \end{aligned}
  \right.
  \end{equation} 
\end{itemize}
Thus $$ U^0= (R^{\infty}(\xi) r_+)(f).$$
So for $f\in H^{2}(\mathbb{R}^n)$, we define $\tilde{R}^{\infty}(\xi)(f)$ by :

\[\tilde{R}^{\infty}(\xi)(f)(x)= \textbf{ }\left\{ 
\begin{array}{l l}
  U^0(x) \quad\text{ if } x\in\Omega_+,\\ 
  \displaystyle \psi(x) \tilde{\alpha}_{0}(x)e^{-M\varphi(x)}\quad\text{if } x\in \Omega_-, & \\ \end{array} \right. \]
where $\tilde{\alpha}_0= \mathcal{E}_{\Omega_-}(P_+U^0)$, $\mathcal{E}_{\Omega_-}$ is the  trace relevement from $\partial\Omega$ to $\Omega_-$ defined on $\displaystyle H^{\frac{3}{2}}(\partial\Omega)$ with values in $H^{2}(\Omega_-)$ and $\psi$ is cut-off function with support in $\overline{\Omega_1}                                                                                                 :=\lbrace x\in \Omega_- : \text{dist}(x,\partial\Omega)\leqslant\delta\rbrace$ and equal to 1 in $\Omega_2 :=\lbrace x\in \Omega_- : \text{dist}(x,\partial\Omega)<\frac{\delta}{2}\rbrace$. Note that $(R^{\infty}(\xi)r_+)(f) =( r_+\tilde{R}^{\infty}(\xi))(f)$.
\begin{remark} Our strategy can be illustrated as follows: we start by determining  $P_-V^1$ and $P_-U^1$ on $\partial\Omega$, which allows us to know $U^1$ thanks to \eqref{213}. Finally, we use (\rm{iv}) and (\rm{viii}) to determine $V^1$.
\end{remark}
\begin{itemize}
\item{ Using (\rm{x}) and the property \eqref{P_+}, one gets }
\[ \textbf{ }\left\{ 
\begin{array}{l l}
  \displaystyle P_-\alpha_{1}(x)= -\frac{1}{2}\Theta_n P_+(T_n-\xi) P_+\alpha_0(x)\quad\text{ in }\partial\Omega,\quad\quad\\
  \displaystyle P_-V^1(x,z)= -\frac{e^{-z}}{2}P_-\Theta_n (T_n-\xi) P_+\alpha_0(x) \quad\text{ in }\partial\Omega\times[0,\delta],\\ 
  \displaystyle P_+V^{1}(x,z)= e^{-z}P_+\alpha_1(x) -ze^{-z}P_+\Theta_n(T_n-\xi) P_+\alpha_0(x)\quad\text{ in }\partial\Omega\times[0,\delta]\quad\quad & \\ \end{array} \right. \]
  
\end{itemize}
Thereby, we obtain that $$V^1(x,z)= e^{-z}\left(P_+\alpha_1(x) -zP_+\Theta_n(T_n-\xi) P_+\alpha_0(x)-\frac{1}{2}P_-\Theta_n (T_n-\xi) P_+\alpha_0(x)\right).$$ 
\begin{itemize}
	\item {Using (\rm{ii}), the definition of $(\tilde{u}^M,\tilde{v}^M)$ and the boundary conditions,  we get that $U^1$ is uniquely} determined by :\end{itemize} 
	\begin{equation}\label{213} 
  \left\{
    \begin{aligned}
      (T_n-\xi)U^1&= 0 \quad\text{ in } \Omega_+,\\ 
\displaystyle P_-U^{1}&= -\frac{1}{2}\Theta_n P_+ (T_n-\xi)P_+\alpha_0 + P_-\Theta_nf\quad\text{on } \partial\Omega ,
    \end{aligned}
  \right.
\end{equation} 
here the term $P_-\Theta_nf$ comes from the contribution of $\omega^M$ (see the proof of Theorem \ref{th1}). The fact that $U^1$ is only determined by \eqref{213} is due to the self adjointness of the Dirac operator $\textbf{H}_\infty$. Indeed, put $G:=\left(U^1-\mathcal{E}_{\Omega_+}(P_-U^{1})\right)\in H^1(\Omega_+)$, then solve the problem \eqref{213} is equivalent to solve  
\begin{equation*} 
  \left\{
    \begin{aligned}
      (T_n-\xi)G&= -(T_n-\xi)\mathcal{E}_{\Omega_+}(P_-U^{1}) \quad\text{ in } \Omega_+,\\ 
\displaystyle P_-G&= 0\quad\text{ on } \partial\Omega_+ .
    \end{aligned}
  \right.
\end{equation*}
Since the self adjointness of $\textbf{H}_\infty$ ensures the existence and uniqueness of $G$, thus $U^1$ exists and it is unique.
Since $P_+V^{1}(x,0)= P_+\alpha_1(x)$, using (\rm{viii}) we get that $P_+\alpha_1 =P_+U^1$ on $\partial\Omega$.
Finally, we define $\tilde{R}^{M}_{1}(\xi)(f)$ by :
\begin{equation}\label{11} 
  \tilde{R}^{M}_{1}(\xi)(f)(x)=\left\{
    \begin{aligned}
      &U^1(x) \quad\text{ if } x\in\Omega_+,\quad\quad\quad\quad\quad\quad\quad\\ 
\displaystyle & \sum_{j=1}^{2}\psi(x) \tilde{\alpha}_{1,j}(x)(M\varphi(x))^{j-1}e^{-M\varphi(x)}\quad\text{if } x\in \Omega_-, \\ 
    \end{aligned}
  \right.
\end{equation}

where $\tilde{\alpha}_{1}=\mathcal{E}_{\Omega_-}(P_+U^1)-\frac{1}{2}P_-\Theta_n (T-\xi)\mathcal{E}_{\Omega_-}(P^+U^0)$ and $\tilde{\alpha}_{2}=-P_+\Theta_n (T-\xi)\mathcal{E}_{\Omega_-}(P^+U^0)$. The cut-off function $\psi$ is used to avoid the problem of non regularity of $\varphi$ far from $\partial\Omega$.\\
The following lemma gathers some properties related to the operator $\tilde{R}^{M}_{1}(\xi)$.

\begin{lemma}\label{lem} The operator $\tilde{R}^{M}_{1}(\xi)$ satisfies:
	\begin{enumerate}
		\item $\tilde{R}^{M}_{1}(\xi)(f)$ only depends on $r_+f$, the restriction of $f$ to $\Omega_+$.
		\item $\tilde{R}^{M}_{1}(\xi)$ does not depend on $M$, furthermore there exists a constant $C$ such that :
		\begin{align}\label{12}
		\displaystyle \forall M,\textbf{ }\forall f\in H^2(\mathbb{R}^n),\quad \displaystyle\Vert\tilde{R}^{M}_{1}(\xi)(f)\Vert_{H^1(\mathbb{R}^n)} \leqslant C\Vert f\Vert_{L^2(\Omega_+)}.
		\end{align}
		\item $\displaystyle P^\infty\tilde{R}^{M}_{1}(\xi)P^\infty=\frac{1}{\lambda^\infty}\mathcal{M}P^\infty$, where $\mathcal{M}$ is a linear map in $E^\infty$. The coeficients of the matrix $\mathcal{M}$ in an orthonormal basis $(f_1,\cdots,f_l)$ of $E^\infty$ are given by the formula :
		$$ m_{kj}=\frac{1}{2}\langle f_j,f_k \rangle_{\partial\Omega}$$
	\end{enumerate}
	
\end{lemma}

\textbf{Proof. } Let $f\in H^2(\mathbb{R}^n)$. We define $U^0$ by \eqref{10}, then theorem \ref{th4} implies that $U^0 \in H^{3}(\Omega_+)$. Thus $U^0$ depends only on $r_+f$. Hence $P_{\pm}\Theta_n( T_n-\xi)(P_+U_{0})_{|\partial\Omega}\in H^{\frac{3}{2}}(\partial\Omega)$ only depends on $r_+f$. Use \eqref{213}, lemma \ref{lemme 1} and the property of $\mathcal{E}_{\Omega_-}$, it provides
\begin{align*}
\Vert\nabla U^1\Vert_{L^2(\Omega_+)} + \Vert U^1\Vert_{L^2(\Omega_+)}\leqslant C \Vert P_- U^1\Vert_{L^2(\Omega_+)} \leqslant C \Vert f \Vert_{L^2(\Omega_+)},\\
\Vert P_+\alpha_1\Vert_{L^2(\partial\Omega)}=\Vert P_+ U^1\Vert_{L^2(\partial\Omega)}\leqslant C\Vert U^1\Vert_{H^1(\Omega_+)}\leqslant C\Vert f\Vert_{L^2(\Omega_+)}.
\end{align*} 
 Moreover, using again the property of $\mathcal{E}_{\Omega_-}$, we get 
\vspace{-1em}
\begin{align*}
\Vert \tilde{R}^{M}_{1}(\xi)(f)\Vert_{L^2(\Omega_-)} &\leqslant \Vert \psi \Vert_{L^\infty(\Omega_-)}\left(\Vert \tilde{\alpha}_1\Vert_{L^2(\Omega_-)}\Vert e^{-M\varphi} \Vert_{L^\infty(\Omega_-)} + \Vert \tilde{\alpha}_2\Vert_{L^2(\Omega_-)}\Vert M\varphi e^{-M\varphi} \Vert_{L^\infty(\Omega_-)}\right)\\
&\leqslant C\left(\Vert \tilde{\alpha}_1\Vert_{H^1(\Omega_-)} + \Vert \tilde{\alpha}_2\Vert_{H^1(\Omega_-)}\right)\\
&\leqslant C\left(\Vert P_+\alpha_1\Vert_{L^2(\partial\Omega)} + \Vert P_-\Theta_n( T_n-\xi)P_+U_{0} \Vert_{L^2(\partial\Omega)}+ \Vert P_+\Theta_n( T_n-\xi)P_+U_{0} \Vert_{L^2(\partial\Omega)}\right)\\
&\leqslant C\Vert f \Vert_{H^2(\Omega_+)}.
\end{align*}
 In addition, we have
\begin{align*}
\Vert \nabla(\tilde{R}^{M}_{1}(\xi)(f))\Vert_{L^2(\Omega_-)} &\leqslant C\Vert f\Vert_{H^2(\Omega_+)}.
\end{align*}
This proves the first and the second assertion.
We now turn to prove the last assertion, i.e $\displaystyle P^\infty\tilde{R}^{M}_{1}(\xi)P^\infty=\frac{1}{\lambda^\infty}\mathcal{M}P^\infty$, where $\mathcal{M}$ is a linear map in $E^\infty$. \\
Let $(f_1,\cdots,f_l)$ be an orthonormal basis of $E^\infty$. Let us compute $P^\infty\tilde{R}^{M}_{1}(\xi)P^\infty(f_j)$. We consider $U^0\in H^3(\Omega_+)$ the solution of :
\[\left\{ 
\begin{array}{l l}
(T_n-\xi)U^0 = f_j \quad\text{ in } \Omega_+,\\ 
\displaystyle P_-U^0 = 0 \quad\text{on }  \partial\Omega. & \\ \end{array} \right. \]
Note that we have $\displaystyle U ^0 = \frac{1}{\lambda^\infty -\xi}f_j$. Let $U^1\in H^2(\Omega_+)$ be the solution of :
\[\left\{ 
\begin{array}{l l}
(T_n-\xi)U^1= 0 \quad\text{ in } \Omega_+,\\ 
\displaystyle P_-U^1=-\frac{1}{2}\Theta_n P_+(T_n-\xi)P_+\alpha_0 + P_-\Theta_n f_j\quad\text{on }  \partial\Omega. & \\ \end{array} \right. \]
We recall that:
$$P^\infty\tilde{R}^{M}_{1}(\xi)P^\infty(f_j) = \sum_{k=1}^{l}\langle U^1,f_k\rangle _{\Omega_+}f_k,  $$
so using the Green's formula we get :
\begin{align*}
 \int_{\Omega_+}U^1f_k &= \frac{1}{\lambda^\infty}\int_{\Omega_+} U^1 T_nf_k\\
 &=\frac{1}{\lambda^\infty}\int_{\Omega_+} T_nU^1 f_k+\frac{i}{\lambda^\infty}\int_{\partial\Omega}(U^1,\boldsymbol{\Lambda_n} f_k)\cdotp n\text{  }d\Gamma\\
  &=\frac{\xi}{\lambda^\infty}\int_{\Omega_+} U^1 f_k+\frac{i}{\lambda^\infty}\int_{\partial\Omega}(P_+U^1 +P_-U^1,\boldsymbol{\Lambda_n} f_k)\cdotp n\text{  }d\Gamma\\
 &=\frac{i}{\lambda^\infty-\xi}\int_{\partial\Omega}(P_+U^1 -\frac{1}{2}\Theta_nP_+(T_n-\xi)P_+\alpha_0 +P_-\Theta_n f_j, \boldsymbol{\Lambda_n}f_k)\cdotp n\text{  }d\Gamma\\
 &=\frac{i}{\lambda^\infty-\xi}\int_{\partial\Omega}(P_+U^1+\frac{1}{2}\Theta_n P_+ f_j, \boldsymbol{\Lambda_n}f_k)\cdotp n\text{  }d\Gamma\\
 &=\frac{i}{\lambda^\infty-\xi}\int_{\partial\Omega}(\Theta_n P_+U^1+\frac{1}{2}P_+ f_j, f_k)\text{  }d\Gamma\\
 &=\frac{1}{2(\lambda^\infty-\xi)}\int_{\partial\Omega}( f_j, f_k)\text{  }d\Gamma,
 \end{align*}
 here we used the lemma \ref{appen}, the fact $ U^1=-\frac{1}{2}\Theta_nP_+(T_n-\xi)P_+\alpha_0 +P_-\Theta_n f_j$ and $(\Theta_n P_+U^1, f_k)=0$ on $\partial\Omega$, where $P_-f_k=0$. Thus  we get :
 $$P^\infty\tilde{R}^{M}_{1}P^\infty(f_j)= \frac{1}{2(\lambda^\infty-\xi)}\sum_{k=1}^{l}\left(\int_{\partial\Omega}( f_j, f_k)\text{  }d\Gamma\right)f_k.$$\qed \\
 Having disposed of this preliminary step, we can now return to the following proposition:

 \begin{proposition}\label{pr23}
The resolvent $\tilde{R}^{M}(\xi)$ admits an asymptotic expansion of the form:
\begin{align}\label{13}
\tilde{R}^{M}(\xi) = e_+R^{\infty}(\xi)r_+ + \frac{1}{M}\tilde{R}^{M}_{1}(\xi) + \frac{1}{M}\tilde{S}^{M}_{2}(\xi).
\end{align}
where $\tilde{R}^{M}_{1}(\xi)$ is defined by \eqref{11} . Furthermore, there is $M_0>0$ and there exists a constant $C$ independent of $\xi \in\mathcal{C}(\lambda^\infty,\eta)$ such that:
\begin{align*}
\sqrt{M}\Vert\tilde{S}^{M}_{2}(\xi)(f)\Vert _{H^1(\Omega_+)}+\Vert\tilde{S}^{M}_{2}(\xi)(f)\Vert _{H^1(\Omega_-)}\leqslant  C\Vert f\Vert_{H^2(\mathbb{R}^n)},\quad\forall M>M_0,\forall f\in H^2(\mathbb{R}^n).
\end{align*}
  \end{proposition}
 
 \textbf{Proof. } Let $f\in H^2(\mathbb{R}^n)$, we set 
 \[\left\{ 
 \begin{array}{l l}
 u^M(x)= U^0(x) + \frac{1}{M}U^1(x) +\frac{1}{M} a^M(x)  \quad\text{ if } x\in\Omega_+,\\ 
 \displaystyle v^M(x)=\psi(x)\tilde{\alpha}_0(x)e^{-M\varphi(x)}+\frac{1}{M}\sum_{j=1}^{2} \psi(x)\tilde{\alpha}_{j}(x)(M\varphi(x))^{j-1}e^{-M\varphi(x)}+ \frac{1}{M}b^M(x)\quad\text{if } x\in\Omega_-. & \\ \end{array} \right. \]
 Thus the remainder terms $a^M$ and $b^M$ satisfy the following system:
 \begin{equation}
 \left\{
 \begin{aligned}
 (T_n-\xi) a^M&=0\text{ in }\text{ }\Omega_+,  \\
 (T_n-\xi) b^M + M\Theta_n b^M &=g^M\text{ in }\text{ }\Omega_-,  \\
 P_+a^M &= P_+b^M \text{ on }\text{ } \partial\Omega,  \\ 
 \end{aligned}
 \right.
 \end{equation} 
 where
 \begin{align*}
 g^M =-e^{-M\varphi}\left( M(T_n\psi )\tilde{\alpha}_0   +(T_n\psi)\sum_{j=1}^{2}\tilde{\alpha}_{j}(M\varphi)^{j-1} -\sum_{j=1}^{2}\psi(T_n-\xi)\tilde{\alpha}_j(M\varphi)^{j-1}-M\Theta_n\sum_{j=1}^{2}\psi\tilde{\alpha}_{j}(M\varphi)^{j-1}\right),
 \end{align*}
so using the fact that $z^je^{-z}$ is bounded and $\varphi(x)>\epsilon>0$ on the support of the derivative of $\psi$,  we get the following estimate :
 $$\Vert g^M\Vert_{L^2(\Omega_-)}\leqslant C\Vert f\Vert_{H^2(\Omega_+)}.$$
Furthermore, if we set $\omega^M = \omega^0 +\frac{1}{M}\omega^1+\frac{1}{M}c^{M}_{1}$ (see the proof of theorem \ref{th1}), we obtain  
\begin{align*}
\Vert P_-b^M\Vert^{2}_{L^2(\partial\Omega)} \geqslant \Vert P_-a^M\Vert^{2}_{L^2(\partial\Omega)} - \frac{C}{M}\Vert f\Vert^{2}_{H^1(\mathbb{R}^n)},
\end{align*}
Hence, applying proposition \ref{prop11}, we get 
 \begin{align*}
\Vert b^M\Vert^{2}_{L^2(\Omega_-)}&\leqslant \frac{C}{M^2}\Vert f\Vert^{2}_{H^2(\mathbb{R}^n)},\\
\Vert b^M\Vert^{2}_{H^1(\Omega_-)}&\leqslant C\Vert f\Vert^{2}_{H^2(\mathbb{R}^n)},\\
\Vert a^M\Vert^{2}_{H^1(\Omega_+)}&\leqslant \frac{C}{M}\Vert f\Vert^{2}_{H^2(\mathbb{R}^n)}.
\end{align*}
Since $ze^{-z}$ is bounded and $\Vert\tilde{\alpha}_0\Vert_{L^2(\Omega_-)}\leqslant C\Vert f\Vert_{L^2(\mathbb{R}^n)}$ we get that  $\Vert M\varphi \tilde{\alpha}_0 e^{-M\varphi}\Vert_{L^2(\Omega_-)}\leqslant C\Vert f\Vert_{L^2(\mathbb{R}^n)}$. So just take 
\[ \tilde{S}^{M}_{2}(\xi)= \left\{ 
\begin{array}{l l}
  \displaystyle a^{M}(x)\text{ if } x\in \Omega_+,\\
  \displaystyle b^{M}(x)+M\varphi(x) \tilde{\alpha}_0(x) e^{-M\varphi(x)}\text{ if } x\in\Omega_-. & \\ \end{array} \right. \]

 This finishes the proof of the proposition.\qed
\setcounter{equation}{0}  
\section{Proof of Theorem \ref{th1} and Theorem \ref{th2}}
This part is devoted to the proof of Theorem \ref{th1} and Theorem \ref{th2}. To that aim, we first prove the following result concerning the projectors
 \begin{proposition}\label{pr31}
 The projector $P^M$ admits an asymptotic expansion of the form 
   $$ P^M =\tilde{P}^\infty  + \frac{1}{M}P^{1,M}+\frac{1}{M}K^{1,M},$$
   where $P^{1,M}$ is a linear operator defined on $H^2(\mathbb{R}^n)$. Furthermore, there exist a constants $C$ and $M_0$ such that for all $M>M_0$ and for all $f\in H^2(\mathbb{R}^n)$  we have :
   \[\quad\left\{ 
 \begin{array}{l l}
\Vert P^{1,M}(f)\Vert _{H^2(\Omega_+)}\leqslant  C\Vert f\Vert_{H^2(\Omega_+)},\\ 
\Vert P^{1,M}(f)\Vert _{L^2(\Omega_-)}\leqslant  C\Vert f\Vert_{H^2(\mathbb{R}^n)},\\
\sqrt{M}\Vert K^{1,M}(f)\Vert _{H^1(\Omega_+)}+\Vert K^{1,M}(f)\Vert _{L^2(\Omega_-)}\leqslant C \Vert f\Vert_{H^2(\mathbb{R}^n)}. & \\ \end{array} \right. \]
Moreover, $r_+P^{1,M}$ is independent on $M$
 \end{proposition}
   \textbf{Proof of Theorem \ref{th1} and proposition \ref{pr31} }
    We show first that the projector $\tilde{P}^M$ admits an asymptotic expansion of the form 
 $$ \tilde{P}^M =e_+P^\infty r_+ +\frac{1}{M}\tilde{P}^{1,M}+\frac{1}{M}\tilde{K}^{1,M},$$
 Inaddition, there exist a constants $C$ and $M_0$ such that for all $M>M_0$ and for all $f\in H^2(\mathbb{R}^n)$  we have  
  \[\quad\left\{ 
 \begin{array}{l l}
\Vert \tilde{P}^{1,M}(f)\Vert _{H^2(\Omega_+)}+\Vert\tilde{P}^{1,M}(f)\Vert _{L^2(\Omega_-)}\leqslant  C\Vert f\Vert_{H^2(\Omega_+)},\\ 
\sqrt{M}\Vert \tilde{K}^{1,M}(f)\Vert _{H^1(\Omega_+)}+\Vert  \tilde{K}^{1,M}(f)\Vert _{L^2(\Omega_-)}\leqslant C \Vert f\Vert_{H^2(\mathbb{R}^n)}. & \\ \end{array} \right. \]
 For this recall that the operators $(\textbf{H}^M - 0\bigoplus\textbf{ H}^{M}_{+})$ and $\textbf{H}^\infty$ are self-adjoint thus the projections $\tilde{P}^{M}$ and $\tilde{P}^{\infty}$ are given by the relations \eqref{pj1}. Hence, according to lemma \ref{lem} and proposition \ref{pr23}, if we set 
 $$ \tilde{P}^{1,M} =\frac{-1}{2i\pi}\int_{\mathcal{C}(\lambda^{\infty},\eta)} \tilde{R}^{M}_{1}(\xi)d\xi,$$
 and
 $$\tilde{K}^{1,M} =\frac{-1}{2i\pi}\int_{\mathcal{C}(\lambda^{\infty},\eta)} \tilde{S}^{M}_{2}(\xi)d\xi,$$
 thus the first assertion and estimates derive directly from lemma \ref{lem} and proposition \ref{pr23}.     
  Now we recall that $\omega^M$ is defined by :
  \begin{equation}\label{wm}
 \left\{
 \begin{aligned}
 \displaystyle  (\textbf{H}^M_+-\xi)\omega^M = f\quad\text{in }\Omega_-  ,&\\
  \displaystyle  P_+\omega^M=0\quad\text{on } \partial\Omega.
 \end{aligned}
 \right.
 \end{equation}
    Then, $\omega^M$ admits an asymptotic expansion of the form 
\begin{align}\label{3.2}  
  \omega^M(x) = \omega^0(x)+ \frac{1}{M}\omega^1(x)+\frac{1}{M}c^M_1(x), 
  \end{align}
 where $\omega^i$ satisfy the following equations :
 
 \[  \left\{ 
\begin{array}{l l}
  \displaystyle  \omega^0 = 0\quad\text{in }\Omega_-  ,&\\
   \displaystyle  (T_n-\xi)\omega^0+ \Theta_n\omega^1 = f\quad\text{in }\Omega_- ,&\\
  \displaystyle  P_+\omega^i=0\quad\text{on } \partial\Omega, \text{ for  }i\in\lbrace0,1\rbrace. & \\ \end{array} \right. \]
  Thus $\omega^0 = 0$ and $\omega^1 = \Theta_nf$.
  Furthermore the remainder term satisfies
   \[  \left\{ 
\begin{array}{l l}
     \displaystyle  (T_n-\xi)\omega^1+ M\Theta_nc^M_1 = 0\quad\text{in }\Omega_- ,&\\
  \displaystyle  P_+c^M_1=0\quad\text{on } \partial\Omega . & \\ \end{array} \right. \]
Here we replaced $\omega^M$ by its asymptotic expansion in the equation \eqref{wm} and we have identified the different powers of $M$. Thus we get $$\omega^M(x) = \frac{1}{M}\Theta_nf(x)+\frac{1}{M}c_1^M, $$
 and we have the estimate
 \begin{align}\label{3.1}
 \Vert c_1^M\Vert_{L^2(\Omega_-)}\leqslant \frac{C}{M}\Vert f\Vert_{H^2(\mathbb{R}^n)}.
 \end{align}
 So, according to proposition \ref{pr23} and the definition of $(\tilde{u}^M,\tilde{v}^M)$, we obtain that
  \begin{align*}
  R^M(\xi) = e_+R^\infty(\xi)r_+ + \frac{1}{M}R^{1,M}(\xi)+\frac{1}{M}S^{2,M},
  \end{align*}
  where the operator $R^M_1(\xi)$ is  defined by: for $f\in H^2(\mathbb{R}^n)$
   \[ R^M_1(\xi)(f)(x)= \left\{ 
\begin{array}{l l}
  \displaystyle U^1(x) \text{ if } x\in \Omega_+,\\
  \displaystyle \psi(x)\tilde{\alpha}_1(x)e^{-M\varphi(x)}+\Theta_nf(x) \text{ if } x\in\Omega_-. & \\ \end{array} \right. \]
  and the operator $S^M_2$ is given by:
  \[ S^M_2(\xi)(f)(x)= \left\{ 
\begin{array}{l l}
  \displaystyle \tilde{S}^M_2(f)(x) \text{ if } x\in \Omega_+,\\
  \displaystyle \tilde{S}^M_2(f)(x)+ c^M_1(x) \text{ if } x\in\Omega_-. & \\ \end{array} \right. \]
  Use the estimate \eqref{3.1}, it provides
$$  \displaystyle \Vert S^M_2(f)\Vert_{H^1(\Omega_+)} +\Vert S^M_2(f)\Vert_{L^2(\Omega_-)}\leqslant C\Vert f\Vert_{H^2(\mathbb{R}^n)}.$$
    This prove Theorem \ref{th1}. If we set 
  $$ K^{1,M} = \frac{-1}{2i\pi}\int_{\mathcal{C}(\lambda^\infty,\eta)}S_2^M(\xi)d\xi.$$
  Then, using the same arguments as in the previous proposition we obtain that
   $$ P^M =\tilde{P}^\infty + \frac{1}{M}P^{1,M}+\frac{1}{M}K^{1,M},$$  
  
  	and we get the estimates looked for. The fact that $r_+P$ is independant on $M$ follows from the first statement of lemma \ref{lem} \qed

  \textbf{Proof of Theorem \ref{th2}.}
  The main idea of the proof is to use Kato's method (see \cite{Kato}). So we introduce the tansformation operator $\mathcal{U}^M$ defined by
  $$ \mathcal{U}^M = I- \tilde{P}^\infty + P^M\tilde{P}^\infty.$$
  This operator maps $\displaystyle E^\infty$ into $\displaystyle \tilde{E}^M$ and leaves the orthogonal of $E^\infty$ invariant. Using the proposition \ref{pr31} we deduce an asymptotic expansion of $\mathcal{U}^M$:
  \begin{align}\label{19}
  \mathcal{U}^M = I +\frac{1}{M}P^{1,M}\tilde{P}^\infty+\frac{1}{M}K^{1,M}\tilde{P}^\infty,
  \end{align}
  We observe that there exists a constant $C$ such that : 
\begin{equation}\label{16}
\displaystyle \forall M>M_0,\quad \forall f\in L^2(\mathbb{R}^n),\quad
 \left\{
 \begin{aligned}
 \displaystyle\left\|  P^{1,M}\tilde{P}^\infty(f)\right\| _{L^2(\mathbb{R}^n)}\leqslant C\left\|  \tilde{P}^\infty(f) \right\| _{H^2(\mathbb{R}^n)}&\leqslant C\left\| f \right\| _{L^2(\mathbb{R}^n)},\\ 
 \displaystyle\left\| K^{1,M}\tilde{P}^\infty(f)\right\| _{L^2(\mathbb{R}^n)}\leqslant  C\Vert \tilde{P}^\infty(f)\Vert_{H^2(\mathbb{R}^n)}&\leqslant C\left\| f \right\| _{L^2(\mathbb{R}^n)}.\quad
 \end{aligned}
 \right.
 \end{equation} 
So by Neumann lemma, we deduce that $\mathcal{U}^M$ is invertible in $\mathcal{L}(L^2)$ and we have
\begin{align*}
(\mathcal{U}^M)^{-1} = I - \frac{1}{M}P^{1,M}\tilde{P}^\infty-\frac{1}{M}K^{1,M}\tilde{P}^\infty+ \sum_{n=2}^{\infty}\frac{(-1)^n}{M^n}\left[ P^{1,M}\tilde{P}^\infty+K^{1,M}\tilde{P}^\infty\right]^n.
\end{align*}
Using estimates \eqref{16} we remark that 
\begin{align*}
\left\| \sum_{n=2}^{\infty}\frac{(-1)^n}{M^n}\left[ P^{1,M}\tilde{P}^\infty+K^{1,M}\tilde{P}^\infty\right]^n \right\|_{\mathcal{L}(L^2)}\leqslant \frac{C}{M^{2}}.
\end{align*}
Hence, we obtain that 
\begin{align*}
\tilde{P}^\infty (\mathcal{U}^M)^{-1} = \tilde{P}^\infty - \frac{1}{M}\tilde{P}^\infty P^{1,M}\tilde{P}^\infty-\frac{1}{M}\tilde{P}^\infty K^{1,M}\tilde{P}^\infty+ \sum_{n=2}^{\infty}\frac{(-1)^n}{M^n}\tilde{P}^\infty\left[ P^{1,M}\tilde{P}^\infty+K^{1,M}\tilde{P}^\infty\right]^n.
\end{align*}
Using the fact that $\tilde{P}^\infty$ is a regularizing operator, we obtain that

\[ \left\{ 
\begin{array}{l l}
\displaystyle\Vert \tilde{P}^\infty P^{1,M}\tilde{P}^\infty\Vert_{\mathcal{L}(L^2)}\leqslant C\Vert P^{1,M}\tilde{P}^\infty\Vert_{\mathcal{L}(L^2 )}\leqslant C,\\
\displaystyle\Vert \tilde{P}^\infty K^{1,M}\tilde{P}^\infty\Vert_{\mathcal{L}(L^2)}\leqslant C\Vert K^{1,M}\tilde{P}^\infty\Vert_{\mathcal{L}(L^2)}\leqslant C.& \\ 
\end{array} \right. \]
So we get that 
\begin{align}\label{17}
\tilde{P}^\infty (\mathcal{U}^M)^{-1} = \tilde{P}^\infty +\frac{1}{M}L^M,
\end{align}
with 
$$ \Vert L^M\Vert_{\mathcal{L}(L^2)}\leqslant C.$$
As a result, the eigenvalues of $\textbf{H}_M$ contained in $B(\lambda^\infty,\eta)$ are the eigenvalues of the operator $\mathcal{Q}^M$ defined on $E^\infty$ by 
$$ \mathcal{Q}^M = \tilde{P}^\infty(\mathcal{U}^M)^{-1}\textbf{H}_MP^M\mathcal{U}^{M}\tilde{P}^\infty.$$
An application of the spectral theorem gives us the relation :
\begin{align}\label{18}
\textbf{H}_MP^M = \frac{-1}{2i\pi}\int_{\mathcal{C}(\lambda^\infty,\eta)}\xi R^M(\xi)d\xi,
\end{align}
so we are reduced to the study of the asymptotic expansion of the operator $\tilde{P}^\infty(\mathcal{U}^M)^{-1}R^MP^M\mathcal{U}^{M}\tilde{P}^\infty$. Using proposition \ref{pr23}, equations \eqref{19} and \eqref{17} we have :
\begin{align*}
\tilde{P}^\infty(\mathcal{U}^M)^{-1}R^MP^M\mathcal{U}^{M}\tilde{P}^\infty  = &\left( \tilde{P}^\infty + \frac{1}{M}L^M\right)\left(\tilde{R}^\infty(\xi) +\frac{1}{M}R^{1,M}(\xi) +\frac{1}{M}S^{2,M}\right)\\
&\times\left( \tilde{P}^\infty +\frac{1}{M}P^{1,M}\tilde{P}^\infty+\frac{1}{M}K^{1,M}\tilde{P}^\infty\right),
\end{align*}
hence
$$\tilde{P}^\infty(\mathcal{U}^M)^{-1}R^MP^M\mathcal{U}^{M}\tilde{P}^\infty = J_1+\ldots+J_6,$$
with
\begin{align*}
J_1 &= \tilde{P}^\infty\tilde{R}^\infty(\xi)\tilde{P}^\infty +\frac{1}{M}\tilde{P}^\infty\tilde{R}^\infty(\xi)P^{1,M}\tilde{P}^\infty +\frac{1}{M}\tilde{P}^\infty\tilde{R}^\infty(\xi)K^{1,M}\tilde{P}^\infty,\quad\quad\quad\quad&\\
\text{ }J_2 &= \frac{1}{M}\tilde{P}^\infty R^{1,M}(\xi)\tilde{P}^\infty + \frac{1}{M^2}\tilde{P}^\infty R^{1,M}(\xi)P^{1,M}\tilde{P}^\infty + \frac{1}{M^2}\tilde{P}^\infty R^{1,M}(\xi)K^{1,M}\tilde{P}^\infty,\\
J_3 &=\frac{1}{M}\tilde{P}^\infty S^{2,M}(\xi)\tilde{P}^\infty + \frac{1}{M^2}\tilde{P}^\infty S^{2,M}(\xi)P^{1,M}\tilde{P}^\infty + \frac{1}{M^2}\tilde{P}^\infty S^{2,M}(\xi)K^{1,M}\tilde{P}^\infty,\\
J_4 &= \frac{1}{M}L^M\tilde{R}^{\infty}(\xi)\tilde{P}^\infty + \frac{1}{M^2}L^M\tilde{R}^{\infty}(\xi)P^{1,M}\tilde{P}^\infty + \frac{1}{M^2}L^M\tilde{R}^{\infty}(\xi)K^{1,M}\tilde{P}^\infty,\quad\quad\\
J_5&= \frac{1}{M}L^M R^{1,M}(\xi)\tilde{P}^\infty + \frac{1}{M^2}L^MR^{1,M}(\xi)P^{1,M}\tilde{P}^\infty + \frac{1}{M^2}L^MR^{1,M}(\xi)K^{1,M}\tilde{P}^\infty,\\
J_6 &= \frac{1}{M}L^M S^{2,M}(\xi)\tilde{P}^\infty + \frac{1}{M^2}L^MS^{2,M}(\xi)P^{1,M}\tilde{P}^\infty + \frac{1}{M^2}L^MS^{2,M}(\xi)K^{1,M}\tilde{P}^\infty.\textbf{  }
\end{align*}
Using lemma \ref{lem}, proposition \ref{pr23}, proposition \ref{pr31}, proposition \ref{pr31} and the fact that $\tilde{P}^\infty$ is a regularizing operator, we estimate the terms $J_1,\ldots,J_6$.
\begin{itemize}
\item $\displaystyle\tilde{P}^\infty \tilde{R}^\infty(\xi)= \frac{1}{\lambda^\infty-\xi}\tilde{P}^\infty$. Thus $\displaystyle J_1= \frac{1}{\lambda-\xi}\tilde{P}^\infty + \frac{1}{\lambda^\infty-\xi}T_1^M$, with $T_1^M =\frac{1}{M}\left(\tilde{P}^\infty P^{1,M}\tilde{P}^\infty +\tilde{P}^\infty K^{1,M}\tilde{P}^\infty\right)$, and we have 
$$ \Vert T_1^M\Vert_{\mathcal{L}(L^2)}\leqslant \frac{C}{M}.$$
\item $\displaystyle J_2=\frac{1}{M} \frac{1}{\lambda-\xi}\mathcal{M}\tilde{P}^\infty + T_2^M$, where $T_2^M =\frac{1}{M^2}\left(\tilde{P}^\infty R^{1,M}(\xi) P^{1,M}\tilde{P}^\infty +\tilde{P}^\infty R^{1,M }(\xi)K^{1,M}\tilde{P}^\infty\right)$ and we have 
$$\Vert T_1^M\Vert_{\mathcal{L}(L^2)}\leqslant\frac{C}{M^2}.$$
\item Using again the fact that $\tilde{P}^\infty$ is a regularizing operator and the proposition \ref{pr31}, we obtain that 
\begin{align*}
\Vert J_3\Vert_{\mathcal{L}(L^2)}\leqslant &\frac{1}{M}\Vert  S^{2,M}\Vert_{\mathcal{L}(L^2)}+\frac{1}{M^2}\Vert \tilde{P}^\infty S^{2,M}\Vert_{\mathcal{L}(L^2)}\\
&\times\left( \Vert P^{1,M}\tilde{P}^\infty\Vert_{\mathcal{L}(L^2)} + \Vert K^{1,M}\tilde{P}^\infty\Vert_{\mathcal{L}(L^2)}\right)\\
&\leqslant \frac{C}{M}.
\end{align*}
\end{itemize}
In the same manner, we can see that $J_4,J_5$ and $J_6$ are bounded in $\mathcal{L}(L^2)$ by $\frac{C}{M}$.\\
In this way, we get 
\begin{align*}
\tilde{P}^\infty(\mathcal{U}^M)^{-1}R^M\mathcal{U}^M\tilde{P}^\infty = \frac{1}{\lambda^\infty - \xi}\tilde{P}^\infty + \frac{1}{M(\lambda^\infty - \xi)}\mathcal{M}\tilde{P}^\infty + \emph{o}(M^{-1}).
\end{align*}
Thus integrating the last asymptotic expansion on $\mathcal{C}(\lambda^\infty,\eta)$, we obtain the following  
\begin{align*}
\tilde{P}^\infty(\mathcal{U}^M)^{-1}\textbf{H}_M\mathcal{U}^M\tilde{P}^\infty = \lambda^\infty\tilde{P}^\infty + \frac{1}{M}\mathcal{M}\tilde{P}^\infty + \emph{o}(M^{-1}).
\end{align*}
Finally, the proof of Theorem \ref{th2} is completed by the classical results on the finite-dimensional perturbation theory (see \cite{Kato},\cite{ReSi4}).
\qed
\setcounter{equation}{0}  
\section{Asymptotic expansion at any order}
In this section, we are looking for an asymptotic expansion of the eigenvalues of $\textbf{H}_M$ at any order. For $M$ sufficiently large, we introduce the operator 
\begin{align}
\label{20}\tilde{\mathcal{U}}^M = (I-W^M)^{-\frac{1}{2}}\left( P^M\tilde{P}^\infty + (I-P^M)(I-\tilde{P}^\infty)\right),
\end{align}
where $W^M := (P^M-\tilde{P}^\infty)^2$ is a self-adjoint operator which tends to zero as $M$ tends to $+\infty$. Moreover $W^M$ commute with $P^M$ and $\tilde{P}^\infty$, thereby
 $\tilde{\mathcal{U}}^M$ maps $\tilde{E}^\infty$ in $E^M$ and $(\tilde{E}^\infty)^{\bot}$ in $(E^M)^{\bot}$. Furthermore, using that $\tilde{P}^\infty\tilde{P}^\infty=\tilde{P}^\infty$ and $P^MP^M=P^M$, we observe that 
 $$\tilde{\mathcal{U}}^{M}\text{ }^{t}\tilde{\mathcal{U}}^M = (I-W^M)^{-\frac{1}{2}}(I-W^M)(I-W^M)^{-\frac{1}{2}}= I,$$
 thus $\tilde{\mathcal{U}}^{M}$ is a unitary, than $\textbf{H}_MP^M$ and $\tilde{\mathcal{Q}}^M = \tilde{P}^\infty(\tilde{\mathcal{U}}^{M})^{-1}\textbf{H}_M\tilde{\mathcal{U}}^{M}\tilde{P}^\infty$ have the same eigenvalues.
 
The next statement follows from \cite[Proposition 5.1]{BrCa02}
\begin{proposition}\label{21}
For  $\tilde{\mathcal{U}}^{M}$ defined above, we have :
$$  \tilde{P}^\infty(\tilde{\mathcal{U}}^{M})^{-1}\textbf{H}_M\tilde{\mathcal{U}}^{M}\tilde{P}^\infty = (I-W^M\tilde{P}^\infty)^{-\frac{1}{2}}\tilde{P}^\infty \textbf{H}_M P^M\tilde{P}^\infty(I-W^M\tilde{P}^\infty)^{-\frac{1}{2}}\tilde{P}^\infty. $$ 
\end{proposition}
This result is important to obtain the asymptotic expansion at any order. Indeed, using the relation \eqref{18} the asymptotic expansion of $\tilde{P}^\infty(\tilde{\mathcal{U}}^{M})^{-1}\textbf{H}_M\tilde{\mathcal{U}}^{M}\tilde{P}^\infty$ will be a consequence of the asymptotic expansion of the operator  $$(I-W^M\tilde{P}^\infty)^{-\frac{1}{2}}\tilde{P}^\infty \textbf{H}_M P^M\tilde{P}^\infty(I-W^M\tilde{P}^\infty)^{-\frac{1}{2}}\tilde{P}^\infty.$$ 
Furthermore, we have that $W^M\tilde{P}^\infty =\tilde{P}^\infty W^M\tilde{P}^\infty= \tilde{P}^\infty- \tilde{P}^\infty P^M\tilde{P}^\infty$, so according to relation \eqref{pj2}, this reduce our discussion to the asymptotic expansion of $P^MR^M(\xi)P^M$. According to Theorem \ref{th4}, we remark that if $f\in H^1(\mathbb{R}^n)$ then $\tilde{P}^\infty(f) \in C^{\infty}(\Omega_+)$ and vanishes on $\Omega_-$, So it is enough to seek an asymptotic expansion of $R^M(\xi)(f)$ with $f\in H^1(\Omega_+)\cap C^{\infty}(\Omega_+)$ and $f=0$ in $\Omega_-$.

\begin{proposition}\label{22}
Let $ N \in \mathbb{N}$ and $f\in L^2(\mathbb{R}^n)$ such that: $f=0$ in $\Omega_-$ and the restriction of $f$ in $\Omega_+$ is in $C^\infty (\Omega_+)$. Then $R^M(\xi)(f)$ admits an asymptotic expansion of the form

\begin{equation*}
\displaystyle R^M(\xi)(f)(x)=
 \left\{
 \begin{aligned}
\displaystyle &\sum_{j=0}^{N}\frac{1}{M^j}U^j(x) +  \frac{1}{M^N}a^M(x) \text{ if } x\in \Omega_+,\\ 
 \displaystyle &\psi(x)\alpha_0(x)e^{-M\varphi} +\sum_{j=1}^{N}\frac{1}{M^j}\left(\sum_{k=1}^{j+1}\psi(x)\alpha_{j,k}(x)(M\varphi)^{k-1}e^{-M\varphi}\right) +  \frac{1}{M^N}b^M(x)  \text{ if } x\in \Omega_-.
 \end{aligned}
 \right.
 \end{equation*}
 Moreover, there exists a constant $C$ independent of $f$ and $M$ such that 
 \begin{equation*}
\displaystyle
 \left\{
 \begin{aligned}
\displaystyle \Vert a^M\Vert_{H^1(\Omega_+)}\leqslant \frac{C}{\sqrt{M}}\Vert f\Vert_{H^{2N}(\Omega_+)},\\
 \displaystyle \Vert b^M\Vert_{H^1(\Omega_-)}\leqslant C\Vert f\Vert_{H^{2N}(\Omega_+)},\\
 \displaystyle \Vert b^M\Vert_{L^2(\Omega_-)}\leqslant \frac{C}{M}\Vert f\Vert_{H^{2N}(\Omega_+)} .
 \end{aligned}
 \right.
 \end{equation*}
\end{proposition}
\textbf{Proof. }
Using the same notations as in section~\ref{sect}, we start by seeking a formal asymptotic expansion of $(u^M,v^M)$ of the form  
\begin{equation*}
 \left\{
 \begin{aligned}
\displaystyle u^M(x) &= U^0(x) +\frac{1}{M}U^1(x) + \frac{1}{M^2}U^2(x)+ \ldots\\ 
 \displaystyle v^M(x) &= V^0(x,M\varphi(x))+\frac{1}{M} V^1(x,M\varphi(x))+\frac{1}{M^2} V^2(x,M\varphi(x))+\ldots
 \end{aligned}
 \right.
 \end{equation*}
 where $(u^M,v^M)$ satisfy the system \eqref{23} with $f=0$ in $\Omega_-$. So replacing $u^M$ and $v^M$ by their asymptotic expansion in \eqref{23} and identifying the different powers terms of $M$, one gets :
 \begin{itemize}
 \item In $\Omega_+$ we have :
 \begin{equation*}
 \left\{
 \begin{aligned}
\displaystyle (T_n-\xi)U^0 &= f \text{ in }\Omega_+,\\ 
 \displaystyle (T_n-\xi)U^1 &= 0 \text{ in }\Omega_+,\\
 \displaystyle &\vdots\\
 \displaystyle (T_n-\xi)U^N &= 0 \text{ in }\Omega_+,
 \end{aligned}
 \right.
 \end{equation*}
 \item In $\Omega_-$ we have :
 \begin{equation*}
 \left\{
 \begin{aligned}
 &\displaystyle V^{0}_{z}=-\mathcal{B}^n_{\Omega_+}V^0\quad\text{ in }\Omega_-,\\
 &\displaystyle (T_{n,x}-\xi)V^j +\Theta_n\mathcal{B}^n_{\Omega_+} V^{j+1}_{z}+\Theta_nV^{j+1}=0\quad\text{ in } \Omega_-, \quad\forall j\in\lbrace 0,\ldots,N\rbrace,\\
 &\displaystyle U^j = V^j \text{ on } \partial\Omega,\quad\forall j\in\lbrace 0,\ldots,N\rbrace.
 \end{aligned}
 \right.
 \end{equation*}
 \end{itemize}
 Using the Duhamel formula we get that :
 \begin{itemize}

 \item{ $\displaystyle V^0(x,z)= \alpha_{0}(x)e^{-z\mathcal{B}^n_{\Omega_+}}$, where $\alpha_0=U^0$ on $\partial\Omega$.}
 \item{ for all $j\in\lbrace 0,\ldots,N\rbrace$, $\displaystyle V^{j+1}(x,z) = e^{-z\mathcal{B}^n_{\Omega_+}}\alpha_{j+1}(x) - \int^{z}_{0}e^{-(z-s)\mathcal{B}^n_{\Omega_+}}\mathcal{B}^n_{\Omega_+}\Theta_n (T_{n,x}-\xi)V^j(x,s)\text{ds}$, where $\alpha_{j+1}=U^{j+1}$ in $\partial\Omega$.}
 \end{itemize}
 Since $R^M(\xi)$ is bounded in $L^2$ this leads us to impose the following condition on $V^j$ : $ P_-V^j = \emph{O}\left(1\right)$. Then we construct $U^j$ and $V^j$ by induction for $j\geqslant 1$. Hence we set $V^j(x,z) = \sum_{k=1}^{j+1}\alpha_{j,k}(x)z^{k-1}e^{-z}$, thus $P_-V^{j}$ satisfy
  \begin{equation*}
 \left\{
 \begin{aligned}
 \displaystyle P_-V^{0}&=0\quad\text{ on }\partial\Omega,\\
 \displaystyle P_-V^1& = P_-\alpha_{1,1}=-\frac{1}{2}P_-\Theta_n(T_n-\xi)\alpha_0\quad\text{and}\quad \alpha_{1,2}=P_+\Theta_n(T_n-\xi)\alpha_0\quad\text{ on } \partial\Omega,\\
 \displaystyle P_-V^2 &= P_-\alpha_{2,1}=-\frac{1}{2}P_-\Theta_n(T_n-\xi)(\frac{1}{4}\alpha_{1,2}-\frac{1}{2}\alpha_{1,1})\quad\text{ on } \partial\Omega,\\
 &\vdots
 \end{aligned}
 \right.
 \end{equation*}
 and $U^j$, solution of 
 \begin{equation}\label{24}
 \left\{
 \begin{aligned}
&\displaystyle (T_n-\xi)U^j =  \text{ in }\Omega_+,\\ 
& \displaystyle P_-U^j = P_-V_j \text{ on }\partial\Omega,
 \end{aligned}
 \right.
 \end{equation}
  Then for all $j\geqslant 1$, $U^j$ is uniquely determined by \eqref{24}. Since $P_+U^j=P_+V^j$ on $\partial\Omega$, this determines the unique $V^j$. 
 For $j=0$ we have
 \begin{equation}
 \left\{
 \begin{aligned}
&\displaystyle (T_n-\xi)U^0 = f \text{ in }\Omega_+,\\ 
& \displaystyle P_-U^0 = 0 \text{ on }\partial\Omega,
 \end{aligned}
 \right.
 \end{equation} 
 thus $U^0 = R^\infty(\xi)r_+$ and $\alpha_{0,1}:=\alpha_0 = P_+U^0$ on $\partial\Omega$.

 \textbf{Properties of the profiles. }Using \eqref{24} and the boundary conditions we remark that if $U^0\in H^{2N +1}(\Omega_+)$ and $\alpha_0\in H^{2N+1}(\Omega_-)$, then $\alpha_{1,k\downharpoonright{\partial\Omega}}\in H^{2N -\frac{1}{2}}(\partial\Omega)$ for all $k\in\lbrace 1,2\rbrace$. Hence using the classical results on elliptic equations we obtain that $U^1\in H^{2N}(\Omega)$ and $\alpha_{1,k}\in H^{2N}(\Omega_-)$. Moreover, there exists a constant $C$ independent of $f$ and $M$ such that:
 $$ \Vert \sum_{k=1}^{2}\alpha_{1,k}\Vert_{H^{2N}(\Omega_-)} +   \Vert U^1\Vert_{H^{2N}(\Omega_+)}\leqslant C \Vert f \Vert_{H^{2N}(\Omega_+)}.$$
Now applying this argument again, we can prove by induction that there exists a constant $C$ independent of $f$ and $M$ such that:
\begin{align}\label{26}
\forall j\leqslant N,\quad \Vert \sum_{k=0}^{j+1}\alpha_{j,k}\Vert_{H^{2(N-j)+1}(\Omega_-)} +   \Vert U^j\Vert_{H^{2(N-j)+1}(\Omega_+)}\leqslant C \Vert f \Vert_{H^{2N}(\Omega_+)}.
\end{align} 
\textbf{Estimation of the remainder terms. }

For any $N\in \mathbb{N}$, we decompose the profiles $u^M$ and $v^M$ in the form 
 
 \begin{equation*}
\displaystyle 
 \left\{
 \begin{aligned}
\displaystyle &u^M=\sum_{j=0}^{N}\frac{1}{M^j}U^j(x) +  \frac{1}{M^N}a^M(x) \text{ if } x\in \Omega_+,\\ 
 \displaystyle &v^M=\psi(x)\alpha_{0,1}(x)e^{-M\varphi(x)}+\sum_{j=1}^{N}\frac{1}{M^j}\psi(x)\sum_{k=1}^{j+1}\alpha_{j,k}(x)(M\varphi(x))^{k-1}e^{-M\varphi(x)} +  \frac{1}{M^N}b^M(x)  \text{ if } x\in \Omega_-.
 \end{aligned}
 \right.
 \end{equation*}
 Moreover, the remainder terms $a^M$ and $b^M$ satisfy the following system :
 \begin{equation}\label{25}
 \displaystyle 
 \left\{
 \begin{aligned}
&\displaystyle (T_n-\xi)a^M=0 \text{ in } \Omega_+,\\ 
 &\displaystyle (T_n-\xi)b^M(x)+M\Theta_nb^M = g^M  \text{ in }  \Omega_-,\\
& \displaystyle a^M = b^M \text{ on } \partial\Omega.
 \end{aligned}
 \right.
 \end{equation}
 where
 \begin{align*}
     g^M =-e^{-M\varphi}\left( M^N(T_n\psi )\alpha_{0,1}   +(T_n\psi)\sum_{j=1}^{N}\sum_{k=1}^{j+1}\frac{(M\varphi)^{k-1}}{M^{j-N}}\alpha_{j,k} -\sum_{k=1}^{N+1}\psi(T_n-\xi)\alpha_{N,k}(M\varphi)^{k-1}-M\Theta_n\sum_{k=1}^{N+1}\psi\alpha_{N,k}(M\varphi)^{k-1}\right),
  \end{align*}
 since $\varphi> \epsilon>0$ on the support of the deriative of $\psi$, using estimate \eqref{26} we obtain that $g^M$ is bounded in $L^2(\Omega_-)$ by $C\Vert f\Vert_{H^{2N}(\Omega_+)}$.

Since we have $P_-a^M=P_-b^M$, then $\displaystyle \Vert P_-b^M\Vert^{2}_{L^2(\partial\Omega)} \geqslant \Vert P_-a^M\Vert^{2}_{L^2(\partial\Omega)} - CM^{-1}\Vert f\Vert^{2}_{L^2(\mathbb{R}^n)}$. Thus an application of proposition \ref{prop11} gives the desired estimates. This finishes the proof. \qed 
 \begin{proposition}\label{28} Let $N\in\mathbb{N}$, then 
  $$ \tilde{P}^\infty R^M(\xi)\tilde{P}^\infty = \frac{1}{\lambda^\infty-\xi}\tilde{P}^\infty + \sum_{j=1}^{N}\frac{1}{M^j}\mathcal{N}^j(\xi)\tilde{P}^\infty + \frac{1}{M^N}\mathcal{K}^{N,M}.$$
  Moreover, if $f\in\tilde{E}^\infty$, then there exists a constant $C$ independent of $f$ and $\xi$ such that
  \begin{equation*}
 \displaystyle 
 \left\{
 \begin{aligned}
&\displaystyle \Vert \mathcal{N}^j(\xi)(f)\Vert_{L^2(\Omega_+)}\leqslant C\Vert f\Vert_{L^2(\Omega_+)},\\ 
 &\displaystyle \Vert \mathcal{K}^{N,M}(\xi)(f)\Vert_{L^2(\Omega_+)}\leqslant \frac{C}{\sqrt{M}}\Vert f\Vert_{L^2(\Omega_+)}.
 \end{aligned}
 \right.
 \end{equation*}
 \end{proposition}
 \textbf{Proof. }We start by observing that if $f\in\tilde{E}^\infty$, then $\tilde{P}^\infty f =f$. Since $\tilde{P}^\infty$ is a regularizing operator, there exists a canstant $C$ such that :
 \begin{align}\label{27}
 \Vert f\Vert_{H^{2N}(\Omega_+)}=\Vert\tilde{P}^\infty f\Vert_{H^{2N}(\Omega_+)}\leqslant C\Vert f\Vert_{L^2(\Omega_+)}.
 \end{align}
 Now we introduce the function $U^0,\ldots,U^N$ defined in the last proposition. For all $j\in\lbrace 1,\ldots,N\rbrace$, we denote by $\mathcal{N}^{j}(\xi)$ the operator defined by
 $$\mathcal{N}^{j}(\xi)(f) :=\tilde{P}^\infty(U^j).$$
 So, using \eqref{26} and \eqref{27} we obtain that there exists a constant $C$ independent of $f$ and $\xi$, such that 
 $$ \Vert \mathcal{ N}^{j}(\xi)(f)\Vert_{L^2(\Omega_+)}\leqslant C\Vert f\Vert_{L^2(\Omega_+)}, \forall j\in\lbrace 1,\ldots,N\rbrace.$$
 Furthermore, using the notations of proposition \ref{22}, if we set $\mathcal{K}^{N,M}(\xi)(f):= \tilde{P}^\infty(a^M)$, then there exist $C$ independent of $f$ and $\xi$, such that 
 $$ \Vert \mathcal{ K}^{N,M}(\xi)(f)\Vert_{L^2(\Omega_+)}\leqslant \frac{C}{\sqrt{M}}\Vert f\Vert_{L^2(\Omega_+)}.$$
 which completes the proof.\qed
 
 Let us now state a corollary of this proposition. For $j\geqslant 1$, we consider the operator $P^j\in\mathcal{L}(\tilde{E}^\infty)$ defined by 
  $$\displaystyle P^j :=\frac{-1}{2i\pi}\int_{\mathcal{C}(\lambda^\infty,\eta)}\mathcal{N}^{j}(\xi)d\xi$$
  then we have 
  \begin{corollary}\label{29}
  For any $N\geqslant 1$
  $$ \tilde{P}^\infty P^M\tilde{P}^\infty =\tilde{P}^\infty + \sum_{j=1}^{N}\frac{1}{M^j}P^j\tilde{P}^\infty + \frac{1}{M^N}P^{N,M},$$
   moreover, there is a constant $C$ such that : $\Vert P^{N,M}\Vert_{\mathcal{L}(L^2(\Omega_+))}\leqslant \frac{C}{\sqrt{M}}$.
  \end{corollary}
  
	 \subsection{  Proof of Theorem \ref{th3}}
  Consider the transformation operator $\tilde{\mathcal{U}}^M$ defined by \eqref{20}, we recall that $\tilde{\mathcal{U}}^M$ maps $\tilde{E}^\infty$ in $E^M$. Thus the eigenvalues of $\textbf{H}_MP^M$ are the same as those of $\tilde{\mathcal{Q}}^M = \tilde{P}^\infty(\tilde{\mathcal{U}}^{M})^{-1}H_M\tilde{\mathcal{U}}^{M}\tilde{P}^\infty$ which is self-adjoint. So, exploiting the finite-dimensional perturbation theory we deduce the asymptotic expansion of the eigenvalues of $\tilde{\mathcal{Q}}^M$ from the one of $\tilde{\mathcal{Q}}^M$ itself. Thus according to  the preoposition \ref{21} we have 
  $$\tilde{\mathcal{Q}}^M = (I-W^M\tilde{P}^\infty)^{-\frac{1}{2}}\tilde{P}^\infty \textbf{H}_MP^M\tilde{P}^\infty(I-W^M\tilde{P}^\infty)^{-\frac{1}{2}}\tilde{P}^\infty.$$
  Using proposition \ref{28}, corollary \ref{29} and the fact that $W^M\tilde{P}^\infty = \tilde{P}^\infty -\tilde{P}^\infty P^M\tilde{P}^\infty$, we deduce that $\mathcal{\tilde{Q}}^M$ admits an asymptotic expansion at any order with respect to $1/M$. Using finite dimensional result for the self-adjoint operator $\mathcal{\tilde{Q}}^M$, we obtain that the eigenvalues of $\textbf{H}_M$ admits an asymptotic expansion at any order. This finishes the proof of Theorem \ref{th3}.\qed 
 \appendix
 \section{ Appendix A}
 In this part, some useful properties of the matrices $ \mathcal{B}^n_{\Omega_+}$ and $P_\pm$ already used are recalled.
 \begin{proposition}\label{pro ape} Let $P_\pm$ as in \eqref{30}. Then
 \begin{itemize}
 \item $\Theta_n P_- =P_+ \Theta_n$ and $\Theta_n P_+ =P_- \Theta_n$ .
 \item $ \mathcal{B}^n_{\Omega_+} P_+= P_+$ and  \textbf{  }$\mathcal{B}^n_{\Omega_+}P_-=-P_-$.
 \item  $\Theta_n\mathcal{B}^n_{\Omega_+} = -i\boldsymbol{\Lambda}_n\cdot\textbf{n}_{\Omega_+}$ and $\lbrace \mathcal{B}^n_{\Omega_+}, \Theta_n\rbrace=0$.
 \end{itemize}
 
 \end{proposition}
 \textbf{Proof. }Use the anticommutator relations $\lbrace \alpha_j,\alpha_k\rbrace= 2\delta_{jk}$ for all $j,k\in\lbrace 1,2,3\rbrace$  and $\lbrace \alpha_j,\beta\rbrace= 0$ ( $\lbrace \sigma_j,\sigma_k\rbrace= 2\delta_{jk}$  for all $j,k\in\lbrace 1,2\rbrace$ and $\lbrace \sigma_j,\sigma_3\rbrace= 0$, $\forall j\in\lbrace 1,2\rbrace$ in the 2D case ).\qed
 \begin{lemma}\label{appen} Let $(f_1,\cdots,f_l)$ be an orthonormal basis of $E^\infty$ and consider $U^0\in H^3(\Omega_+)$ the solution of the problem:
\[\left\{ 
\begin{array}{l l}
(T_n-\xi)U^0 = f_j \quad\text{ in } \Omega_+,\\ 
\displaystyle P_-U^0 = 0 \quad\text{on }  \partial\Omega, & \\ \end{array} \right. \]
such that $P_+\alpha_0 = P_+U^0$. Then, the following equality hold in $\partial\Omega$
$$P_+(T_n-\xi)P_+\alpha_0 = P_+f_j.$$ 
 \end{lemma}
 \textbf{Proof.} We prove this result in the case of dimension 2 and 3 separately. All computations are done on $\partial\Omega$.
\begin{itemize}
\item \textbf{2D case}: Since $\alpha_0= P_+\alpha_0$ then $\alpha_0= \begin{pmatrix}
 1\\
 a_2
 \end{pmatrix}
 \alpha_{0,1}$, where $\alpha_{0,1}$ is a scalar function defined on the boundary.
 \end{itemize}
  Using the definition of tubular coordinates (see notation \ref{note6}) and the fact that $\textbf{n}_{\Omega_+}(s)=\left(\cos\theta(s),\sin\theta(s)\right)$, we get
 $$ \begin{pmatrix}
 \partial_t\\
 (1+t\theta^{\prime})^{-1}\partial_s
\end{pmatrix} =\begin{pmatrix}
 \cos\theta & -\sin\theta\\
 \sin\theta & \cos\theta
\end{pmatrix} \begin{pmatrix}
 \partial_1\\
 \partial_2
\end{pmatrix} .$$ 
A direct computation gives 
\begin{align*}
P_+(T_2-\xi)\alpha_0= \frac{1}{2}\begin{pmatrix}
1 \\
a_2
\end{pmatrix}\left( D_2^\ast a_2+a_2^\ast D_2) -2\xi\right)\alpha_{0,2}=\frac{1}{2}\begin{pmatrix}
1 \\
a_2
\end{pmatrix}\left( -2i\partial_s +\theta^{\prime}-2\xi\right)\alpha_{0,2},\\
P_-(T_2-\xi)\alpha_0= \frac{1}{2}\begin{pmatrix}
1 \\
-a_2
\end{pmatrix}\left( D_2^\ast a_2-a_2^\ast D_2\right)\alpha_{0,2}=\frac{1}{2}\begin{pmatrix}
1 \\
-a_2
\end{pmatrix}\left( 2\partial_t +\theta^{\prime}\right)\alpha_{0,2}.
\end{align*}
We note that the operator $P_+(T_2-\xi)P_+$ contains only the $\partial_s$ derivative and $P_-(T_2-\xi)P_-$ contains only the $\partial_t$ derivative. We will show that the same result hold true for the 3D case. 
\begin{itemize}
\item \textbf{3D case}: as in the 2D case we have  $$\alpha_0 = \begin{pmatrix} 
 \alpha_{0,2}\\
 a_3\alpha_{0,2}
 \end{pmatrix}
 := \begin{pmatrix} 
 I_2\\
 a_3
 \end{pmatrix}
 \alpha_{0,2},$$
 \end{itemize} where $I_2$ is the $2\times 2$ identity and $\alpha_{0,2}$ is $1\times2$ vector function defined on the boundary. Let $U$ ( resp. $V$) be an open subset of $\mathbb{R}^2$ ( resp. $\mathbb{R}^3$) and $X: U\rightarrow V$ be a local regular parametrization  of $V\cap \partial\Omega$.
 
 For $s=(s_1,s_2)\in U$, we denote by $g(0)$ the matrix of the first fundamental form on $\partial\Omega$ with entries :
 $$ g(0)(s_1,s_2) = \left(\frac{\partial X}{\partial s_i}\cdot \frac{\partial X}{\partial s_j}\right)(s_1,s_2).$$
 We denote $g^{ij}(0)$ the coefficents of the matrix $g(0)^{-1}$. We know that $g(0)$ and $g(0)^{-1}$ are symmetric.
 We denote by $\nabla_{\partial\Omega}=g^{-1}\begin{pmatrix}
 \partial_{s_1}\\
 \partial_{s_2}
 \end{pmatrix}:=\begin{pmatrix}
 G_1(s)\\
 G_2(s)
 \end{pmatrix}$, $\textbf{n}_{\Omega_+}(s)=(n_1(s),n_2(s),n_3(s))^t$, $\frac{\partial X}{\partial s_1}(s)=(X^1_1(s),X^2_1(s),X^3_1(s))$ and $\frac{\partial X}{\partial s_2}(s)=(X^1_2(s),X^2_2(s),X^3_2(s))$. If we denote by $\boldsymbol{\tilde{\sigma}}:=(\sigma_1,\sigma_2,\sigma_3)$, then the operators $D $ and $a_3$ are written as follows 
 $$D_3=-i\boldsymbol{\tilde{\sigma}}\cdot\nabla=\begin{pmatrix}
 -i\partial_3 & -i\partial_1 + \partial_2\\
 -i\partial_1 - \partial_2 & i\partial_3
\end{pmatrix},\quad\quad a_3=i\boldsymbol{\tilde{\sigma}}\cdot\textbf{n}_{\Omega_+} = \begin{pmatrix}
i n_3& i n_1 + n_2 \\
i n_1 - n_2 & -i n_3
\end{pmatrix},$$  
 furthermore, we have
 $$\begin{pmatrix}
 \partial_1\\
 \partial_2\\
 \partial_3
 \end{pmatrix}=\begin{pmatrix}
 n_1 & X^1_1 &  X^1_2\\
 n_2 & X^2_1 &  X^2_2\\
 n_3 & X^3_1 &  X^3_2
 \end{pmatrix}\begin{pmatrix}
 \partial_t\\
 \nabla_{\partial\Omega}
 \end{pmatrix}.$$
 As in the 2D case, a simple computation gives us 
 \begin{align*}
P_+(T_3-\xi)\alpha_0= \frac{1}{2}\begin{pmatrix}
1 \\
a_3
\end{pmatrix}\left( D_3^\ast a_3+a_3^\ast D_3) -2\xi\right)\alpha_{0,3}=\frac{1}{2}\begin{pmatrix}
1 \\
a_3
\end{pmatrix}\left( L(s)\cdot\nabla_{\partial\Omega}  + D_3^\ast( a_3)-2\xi\right)\alpha_{0,2},\\
P_-(T_3-\xi)\alpha_0= \frac{1}{2}\begin{pmatrix}
1 \\
-a_3
\end{pmatrix}\left( D_3^\ast a_3-a_3^\ast D_3\right)\alpha_{0,3}=\frac{1}{2}\begin{pmatrix}
1 \\
-a_3
\end{pmatrix}\left(2\partial_t I_2+ D_3^\ast (a_3)\right)\alpha_{0,2},
\end{align*}
here we used the properties of Pauli matrices. The operator $L(s)\cdot \nabla_{\partial\Omega}$ that depends only on $\partial_{s }$ is given by  
\begin{align*}
-2i\begin{pmatrix}
(n_1X^2_1-n_2X^1_1)G_1(s) + (n_1X^2_2-n_2X^1_2)G_2(s) & -i\sum^{2}_{j=1}\left[ n_3(X^1_j-iX^2_j)-X^3_j(n_1-in_2)\right]G_j(s) \\
i\sum^{2}_{j=1}\left[ n_3(X^1_j+iX^2_j)-X^3_j(n_1+in_2)\right]G_j(s)  & -n_1X^2_1-n_2X^1_1)G_1(s) - (n_1X^2_2-n_2X^1_2)G_2(s)
\end{pmatrix}.
\end{align*}
Now we go back to the proof of our statement for the 2D case. We remark that any function $U$ may be written as:
\begin{align*}
 U = P_+U +P_-U=\begin{pmatrix}
1\\
a_2
\end{pmatrix}u_+ +\begin{pmatrix}
1\\
-a_2
\end{pmatrix}u_-. 
\end{align*}
Sincee $T_2U =T_2\begin{pmatrix}
1\\
a_2
\end{pmatrix}u_+ + T_2\begin{pmatrix}
1\\
-a_2
\end{pmatrix}u_-$, So for $U$ such that $P_-U=0$ we get that $u_-=0$ and the following holds
\begin{align*}
T_2U&= P_+T_2\begin{pmatrix}
1\\
a_2
\end{pmatrix}u_+ + P_-T_2\begin{pmatrix}
1\\
a_2
\end{pmatrix}u_+\\
&= \begin{pmatrix}
1\\
a_2
\end{pmatrix}\left( -i(1 +t\theta^{\prime})^{-1}\partial_s +\frac{1}{2}\theta^{\prime}(1+t\theta^{\prime})^{-1}\right)u_+ + \begin{pmatrix}
1\\
-a_2
\end{pmatrix}\left( \partial_t +\frac{1}{2}\theta^{\prime}(1+t\theta^{\prime})^{-1}\right)u_+.
\end{align*}
In particular, Every eigenfunction $f_j$ of $\textbf{H}_\infty$ may be written as $\begin{pmatrix}
1\\
a_2
\end{pmatrix}f^+_j$ and satisfies:

\begin{align*}
T_2\begin{pmatrix}
1\\
a_2
\end{pmatrix}f^+_j = \lambda^\infty \begin{pmatrix}
1\\
a_2
\end{pmatrix}f^+_j &\Longleftrightarrow 
 \left\{
 \begin{aligned}
\left(-i(1 +t\theta^{\prime})^{-1}\partial_s + \partial_t+\theta^{\prime}(1+t\theta^{\prime})^{-1}\right) f^+_j &= \lambda^\infty f^+_j  ,\\
\displaystyle \left(-i(1 +t\theta^{\prime})^{-1}\partial_s -\partial_t\right) f^+_j&= \lambda^\infty f^+_j.
 \end{aligned}
 \right. \\ 
 &\Longleftrightarrow 
 \left\{
 \begin{aligned}
\left(\partial_t+\frac{1}{2}\theta^{\prime}(1+t\theta^{\prime})^{-1}\right) f^+_j &= 0,\\
\displaystyle \left(-i(1 +t\theta^{\prime})^{-1}\partial_s +\frac{1}{2}\theta^{\prime}(1+t\theta^{\prime})^{-1}\right) f^+_j&= \lambda^\infty f^+_j.
 \end{aligned}
 \right.  
\end{align*}

Finally, using the fact that $U^0=\frac{1}{\lambda^\infty-\xi}f_j$ and take $t = 0$ in the last equality one gets the desired result. 
This finishes the proof of the 2D case, the detailed verification of the 3D case being left to the reader.\qed

   \textbf{      Acknowledgement.} The author wishes to express his gratitude to his thesis advisor Vincent Bruneau for suggesting the
problem and many stimulating conversations, without which this article would not have been completed.



\emph{Universit\' e de Bordeaux, IMB, UMR 5251, 33405 Talence Cedex (FRANCE)}\\
 \emph{and }\\
 \emph{Departamento de Matem\' aticas, Universidad del Pa\' is Vasco/Euskal Herriko Unibertsitatea, UPV-EHU, 48080 Bilbao (SPAIN)}\\
\emph{E-mail address}: allal.bader@yahoo.fr


\begin{thebibliography}{99}
\bibitem{A} R.A. Adams, {\em Sobolev Spaces}, Pure and Applied Mathematics
{\bf 65}, Academic Press.
\bibitem{AB} A. R. Akhmerov and C. W. J. Beenakker, {\em Boundary conditions for Dirac fermions on a terminated honeycomb lattice}, Phys. Rev. B 77 (2008), 085423.
\bibitem{ATR1} N. Arrizabalaga, L. Le Treust and N. Raymond, {Extension operator for the MIT bag model}, HAL e-prints hal-01540149, 2017.
\bibitem{ATR} N. Arrizabalaga, L. Le Treust and N. Raymond, {\em The MIT Bag Model as an Infinite Mass Limit}, Preprint arXiv: 1808.09746v1, August 2018.
\bibitem{RC} J.-M. Barbaroux, H. Cornean, L. Le trust and E. Stockmeyer, {\em Resolvent convergence to Dirac operators on planar domains}, Priprint arXiv:1810.02957, 2018.
\bibitem{BFSV} R.D. Benguria, S. Fournais, E. Stockmeyer, H. Van Den Bosch, {\em Selfadjointness of two-dimensional Dirac operators on domains}. Ann. Henri Poincar\' e 18(4), 1371-1383, 2017.
\bibitem{BM} M. V. Berry and R. J. Mondragon, {\em Neutrino billiards: time-reversal symmetry-breaking without magnetic fields}. Proc. Roy. Soc. London Ser. A, 412(1842):53-74, 1987.
\bibitem{Bo} P. Bogolioubov. {\em Sur un mod\' ele \' a quarks quasi-ind\' ependants}. Annales. I.H.P., section A, 8:163-189, 1968.
\bibitem{BrCa02} V. Bruneau, G. Carbou, {\em Spectral asymptotic in the large coupling limit}, Asymp. Anal, {\bf 29} (2)  (2002) 91-113.
\bibitem{CGPNG} A. H. Castro Neto, F. Guinea, N. M. R. Peres, K. S. Novoselov, and A. K. Geim, {\em The electronic properties of graphene}, Rev. Mod. Phys. 81 (2009), 109–162.
\bibitem{CJT}  A. Chodos, R. L. Jaffe, K. Johnson and C. B Thorn . {\em Baryon structure in the bag theory}. Phys. Rev. D 10 (8), 2599-2604 (1974).
\bibitem{CJTW} A. Chodos, R. L. Jaffe, K. Johnson, C. B Thorn and V. F Weisskopf. {\em New extended model of hadrons}. Phys. Rev. D 9 (12), 3471-3495 (1974).
\bibitem{Gil} D. Gilbarg and N.S. Trudinger, {\em Elliptic partial differential equations of second order}, volume 224 of Grundlehren der Mathmatischen Wissenschaften [Fundamental Principles of Mathmatical Sciences]. Springer-Verlag, Berlin, Second edition, 1983.
\bibitem{Kato} T. Kato, { \em Perturbation Theory for Linear Operators}, Springer Verlag, Berlin/Heidelberg/New York, 1966. 
\bibitem{MTY} Y. Miake, H. Tetsuo and K. Yagi.{\em Quark-gluon plasma: From big bang to little bang}, Cambridge University Press, 2005.
\bibitem{OV} T. Ourmi\`eres-Bonafos, L. Vega, {\em  A strategy for self-adjointness of Dirac operators: applications to the MIT bag model and delta-shell interactions}. Publ. Mat, 2018.
\bibitem{STSV} E. Stockmeyer and S. Vugalter, {\em Infinit mass boundary conditions for Dirac operators}, Preprint arXiv 1603.09657. Accepted for publication in Journal of Spectral Theory, 2017.
\bibitem{ReSi4} M. Reed, B. Simon, {\em Methods of Modern Mathematical Physics} vol. IV, Analysis of Operators. Academic Press, New York, 1978.
\bibitem{Tha} B. Thaller, {\em The Dirac equation}, Text and Monographs in Physics, Springer-Verlag, Berlin, 1992.




\end{thebibliography}
\end{document}